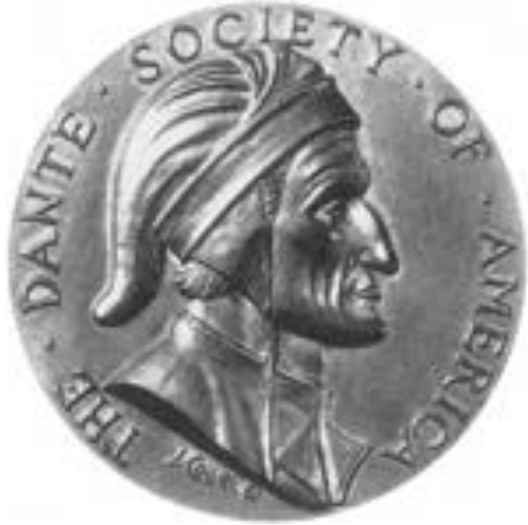

# "The three *giri* of *Paradiso* XXXIII"
Arielle Saiber and Aba Mbirika
*Dante Studies* **131** (2013): 237-272

**ABSTRACT**

Our paper offers an analysis of how Dante describes the *tre giri* ("three rings") of the Holy Trinity in *Paradiso* 33 of the *Divine Comedy*. We point to the myriad possibilities Dante may have been envisioning when he describes his vision of God at this final stage in his journey. Saiber focuses on the features of shape, motion, size, color, and orientation that Dante details in describing the Trinity. Mbirika uses mathematical tools from topology (specifically, knot theory) and combinatorics to analyze all the possible configurations that have a specific layout of three intertwining circles which we find particularly compelling given Dante's description of the Trinity: the round figures arranged in a triangular format with rotational and reflective symmetry. Of the many possible link patterns, we isolate two particularly suggestive arrangements for the *giri*: the Brunnian link and the (3,3)-torus link. These two patterns lend themselves readily to a Trinitarian model.

# The Three *Giri* of *Paradiso* 33

ARIELLE SAIBER AND ABA MBIRIKA

As he gazes deeply into the point of light emanating from the Empyrean at the end of his journey (*Par.* 33.112–14), Dante the Pilgrim feels as if it were changing, although really he knows that it is *he* who is changing. This Divine Light, he lets us know straight away—even with the changes and diverse features he is seeing, and even with its three components he is about to see—is always *una sola parvenza* (113). What he saw as a point of light earlier in the Primum Mobile and on his entrance to the Empyrean he has seen take numerous forms; rather, he has been able to see the true form more and more perfectly and completely. Now, in the *alto lume* (116) of the poem's final tercets, he sees the Divinity as the Trinity.

> Ne la profonda e chiara sussistenza
> de l'alto lume parvermi tre giri
> di tre colori e d'una contenenza;
> e l'un da l'altro come iri da iri
> parea reflesso, e 'l terzo parea foco
> che quinci e quindi igualmente si spiri.
> (*Par.* 33.115–20)[1]

It (or "they," as we will consider below: *parvemi* or *parvermi*?) appear in the form of (or looked like to him) three *giri* (circles/discs, spheres/balls, tori, cylinders, spirals, ellipsoids, or other round things?), *di tre colori* (each a different color, or each containing three colors?), and *d'una contenenza* (all the same size, all occupying a single space, all contained within a single space, and/or all of the same substance?). One of the *giri* appears as a reflection of another, *come iri da iri* (resembling a double rainbow, or merely an analogy for generation?), and the third *giro* appears as a flame (*parea foco*) breathed or breathed forth in equal measure (*si spiri*) by the first two (*quinci e quindi*).

While there is little doubt which person of the Trinity Dante associated with which *giro* (the first *giro*, the Father, generates the second *giro*, the Son, which reflects the first, while both breathe forth the Holy Spirit, the fiery one), the enigmas embedded within the details of these verses are many. How are these three *giri*—very likely circles, given that Dante calls the one associated with the Son a *circulazion* (127) and a *cerchio* (138) (but what "circle" means is not entirely clear)—situated with respect to one another from the Pilgrim's viewpoint? Are they arranged along a horizontal axis or a vertical axis? Are they in a triangular format? Or are they in another arrangement entirely? Are they proximate to one another but not touching, or are they tangent, intersecting, fully overlapping, or interwoven? Are they moving: rotating around their centers, or around each other, or around Dante? Is Dante orbiting them? And what is their chromatic nature? Do they each exhibit a single color, or do they each display the same spectrum of colors?

Dante's depiction of the "three-in-oneness" of the Holy Trinity, while alluded to earlier in the *Commedia* and elsewhere in his works,[2] is one of the very last imagistic puzzles used to express God's ineffability in the poem, and it is also certainly one of the most perplexing. Virgil's pronouncement in *Purgatorio* 3 that "matto è chi spera che nostra ragione / possa trascorrer la infinita via / che tiene una sustanza in tre persone" (34–36) unequivocally alerts us to the impossibility of understanding God's triune nature with reason alone. A number of commentators have, in fact, emphasized the futility of trying to visualize the Trinity Dante describes.[3] The Pilgrim, however, seems to have a theophany and to "see" and understand the Trinity, or at least its similitude, by means of a Divine flash (*fulgore*) that passes through his mind (140–41). The poet made an effort to furnish this vision—abstract though it was and the less-than-*poco* of it that he could recall/articulate (see 106–8 and 121–23)—with specifics. As Mario Fubini observes, the more sublime a topic or object, the more precise the language Dante sought to describe it;[4] or as Carlo Grabher remarks, Dante succeeds in obtaining "tutta la concretezza possibile nella maggiore levità possibile."[5] The poet could have merely echoed his earlier verses describing the Trinity, especially the one in response to St. Peter's inquisition on faith: "e credo in tre persone etterne, e queste / credo una essenza sì una e sì trina, / che soffera congiunto

'sono' ed 'este'" (*Par.* 24.139–41), but he did not. If then it was important enough to the poet to depict the Trinity in its particulars—even if these particulars could be considered "accidents" of the True Substance packaged just for Dante's mortal eyes—then it is important for us to try "to conceive [this] concrete image," as Charles Singleton remarks in his commentary on the *giri*.[6]

Fortunately, not all scholars have abdicated in the face of this crux. They have recognized that by articulating his vision of the Trinity as he did, Dante joined the ranks of myriad medieval theologians, mystics, and artists who devoted themselves to visualizing this paradox, and of thinkers, such as Thomas Aquinas and Bonaventure, and later Eckhart, who explored the means by which the mind can approach knowledge of God *through* reason, even if full understanding cannot be granted without grace (the *fulgore*).[7] Like the difficult *nodo* binding the mystery of the universe's shape that the Pilgrim struggled to untie while in the Primum Mobile (*Par.* 28.58), the Trinity is quite literally a knot, or rather *the* knot, that beckons Dante, and us, into contemplation precisely because it is beyond reason's grasp.

The paradox of three-in-one remains Christianity's greatest mystery. This essay does not engage the additional layer of paradox woven into verses 33.127–33 as to how *la nostra effige*" (131), *l'imago* (138), "fits into" the Son's *giro*; nor does it discuss the geometer's failure to square the circle.[8] It focuses instead on the details of shape, motion, size, and color with which Dante describes the Trinity in order to consider possible geometric configurations of the *tre giri* that may have inspired the poet. In the first five sections we look at the numerous ways that verses 115–120 have been and could be read. In the final section, with tools from the mathematical fields of topology (specifically, knot theory) and combinatorics, we analyze one particularly suggestive arrangement for the *giri*: that of three intertwined circles in a triangular format. Of the many permutations of this figure, we isolate two variations—a Brunnian link commonly called the Borromean rings and a (3,3)-torus link—to show how they more than any other possible arrangement offer unique mathematical, aesthetic, and metaphoric properties that resonate with many of the qualities of the Trinity Dante allusively described in *Paradiso* 33. We propose these as a possible configurations, rich with mystery in themselves, out of a number of Trinitarian models that Dante knew and contemplated.

## Shape

Many scholars have discussed at length the challenges Dante faced as he wrote about the timeless and spaceless simultaneity of the Empyrean through human language, which, as Teodolinda Barolini notes, is time-based and linear.[9] And we agree with Guy Raffa that "the ineffability of paradise means, paradoxically, that Dante's words *become* the poem's content as well as the poet's means of expression" in his incarnational poetry.[10] It is with awareness of these peculiarities, and the recognition that the poet repeatedly makes it clear that his words can only approximate what the Pilgrim saw and experienced and whose memory can recall, that we approach our first conundrum, that is, the meaning of word *giri*. While the term *giro* unquestionably signifies "roundness," it is not a given that Dante intended the *giri* to be "circles" (even with the seemly clear reference to the *cerchio* of the Son in verse 138), as many commentators and translators have construed; nor can it be certain that he intended the *giri* to mean "spheres," as several scholars have postulated as an alternative to "circles."[11] Interestingly, most commentators writing in Latin and Italian up through the late nineteenth century refrained from even discussing the form of the *giri*,[12] and only a few scholars in the last seven centuries have problematized the term to consider its multiple, potential shapes.[13]

Throughout the *Commedia* the substantive *giro* (also appearing as *giri*, *giron[e]*, and *gironi*) has indicated various kinds of round things: a circular region, an object that encircles (a belt, wreath, crown, ring), a twisting/turning/revolving/wheeling, a celestial sphere, a horizon, and a passing of time (that is, a turning of the heavens).[14] And curiously, in the *Convivio*, when speaking about the heaven of Jupiter and its associated art, geometry, Dante wrote that when he says *cerchio* he implies "largamente ogni ritondo, o corpo o superficie" (2.13).[15] If the word *cerchio* could signify a circle, a sphere, or anything round in Dante's mind, then the flexibility of the already polyvalent *giro* jumps further into relief. The three *giri* could be circles, but they could also be discs, spheres/balls, or even tori, cylinders, spirals, ellipsoids, or other round things. What is more, the fact that Dante states that the *giri* "appeared" to the Pilgrim (or "appeared to be"— the verb *parere* could indicate either) may imply that they manifested themselves as they did specifically for his human faculty of vision. As John Carroll says, "the mere fact that Dante uses a figure to convey what he saw of the Trinity is possibly his way of saying that he did not see the Divine Essence, since the very idea of seeing the Essence is face to face vision, without the need of any similitude."[16] What the *giri* actually *are* may, in fact, be entirely other. This sort of double-appearing has happened throughout the *Commedia*, most pertinently to our current discussion in *Paradiso* 28, when the Pilgrim looks up at the God-point from the primum mobile and sees nine concentric *cerchi* of light encircling it (which turn out to be bands of angels); and in *Paradiso* 30, after his eyes drink of the water-waves of light—*luce viva* (49), *rivera* (61), *fiume* (76), *onda* (86), *gronda* (88)—emanating from the God-point, he seems to see this flowing light change its shape from a length/straightness, to a roundness: *di sua lunghezza divenuta tonda* (90).[17]

The "appearing" and "appearing as" duality is even further complicated in the description of the *giri* by the fact that early manuscript copies of the *Commedia* employ two variations of the verb *parere* for verse 116: the plural *parvermi* and the singular *parvemi*. The former emphasizes the "threeness" of the Trinity; the latter, its "oneness." For the Edizione Nazionale of the *Commedia*, Petrocchi opts for *parvermi*, and most subsequent editors follow suit. Either option works conceptually and metrically within the poem, although there is a case to be made in favor of the singular form, as it forces the reader to recall the oneness of the Trinity, which the Pilgrim is seeing manifest as three roundnesses of three colors. With the adoption of the Nicene Creed of 325—which rejected the Arian belief in the subordination of God the Son to God the Father, as the Son proceeds from the Father—, the three persons of the Trinity are made *consubstantiales et coaequales*, of the same substance/essence and coequal. A few verses after describing the appearance of (or appearing to the Pilgrim of) the *tre giri*, Dante refers to the *giro* of the Son as a *circulazion* (127), a circling or circulation. He states two verses later that the *circulazion* was *circunspetta* by his eyes (129), which could be interpreted to mean "carefully looked around or along with," that is, that his eyes moved in a circular path to see the *giri* because they were moving. What begins to become apparent is that the actual shape of the *giri* could be one or more of these "roundnesses" simultaneously, especially if the *giri* were in motion. A ring (circle) spinning on its edge, for example, could look like a sphere revolving about its center, like a wheel or a disc. And depending on the figures' configuration, other "roundnesses" could emerge to the eye. For now, what we can surmise is that the three *giri* "appeared as" (and "appeared to be") circulations. Evidently, the poet did not want to commit to naming a precise geometric form. Perhaps he wished to retain the air of mystery inherent in the Trinity; perhaps he did not want to risk saying something about the Trinity that could be deemed heretical; perhaps he wished to indicate that the Pilgrim was not able to individuate any precise circular form; or perhaps the poet wanted his readers to realize that the form the Pilgrim saw was not the Trinity's true form but only one constructed for human eyes. Whatever the reason for not specifying the *giri*'s form, Dante chose to depict the three persons of the Trinity as geometric "roundnesses" rather than anthropomorphically or with other common iconographic conventions. Such a decision seems a natural extension of his belief that the circle was the most perfect of forms (*Conv.* 2.13), that geometry was "sanza macula d'errore" (*Conv.* 2.13), and that God was "Colui che volse il sesto" (*Par.* 19.40) in order to create the universe. The final canticle of the *Commedia*, and especially the final canto, is a monument to the use of geometric imagery to describe the ineffable. Dante's use of the circle and the sphere as metaphors for the incommensurable, paradoxical, eternal nature of God is part of a long tradition in Western thought.[18]

In his depiction of the Trinity as three "roundnesses," Dante was following an established motif for his time. Theologians such as Pseudo-Dionysius, Hugh and Richard of St. Victor, Alan of Lille, and Thomas discussed the infinity and eternity of the circle/sphere as akin to that of the Trinity. And as early as Saint Augustine's *De Trinitate* (9.4.7) and Petrus Alfonsi's *Dialogi contra Iudaeos* (Tit. 6), the Trinity has been represented as three intersecting rings (fig. 1). One of the most explicit depictions of the Trinity as three circles, however, is that of a twelfth-century diagram thought by many scholars to be conceived by the Franciscan abbot, Joachim of Fiore (fig. 2). The image, table 11 in the *Liber figurarum*,[19] consists of three rings aligned horizontally such that any pair of rings are linked. The first ring is painted green, the second blue, and the third red, with the Latin word for Father, Son, and Holy Spirit written at the top of each ring, respectively. Joachim explains in his *Expositio in Apocalypsim* that the green ring represents the Father; the blue, the Son; and the red, the Holy Spirit (see fol. 36v and 101r).

Leone Tondelli, Marjorie Reeves, and Beatrice Hirsch-Reich have made a strong case for Dante having seen this figure—especially given his studies with Joachimites Pier di Giovanni Olivi and Ubertino da Casale at Santa Croce.[20] Like many scholars of both Joachim and Dante, they believe that Dante's *giri* were partially modeled on Joachim's image—"partially" because Joachim's diagram (along with a number of other aspects of his doctrine and propositions for church reform) was considered erroneous by the Fourth Lateran Council of 1215. Attributing a different color to each circle as Joachim did amounts to attributing a distinction in substance between the three persons. Yet even with the problematic reputation of Joachim in the thirteenth and early fourteenth centuries, Dante places the abbot among the saints and theologians in the Sphere of the Sun in *Paradiso* 12.140–41, at the side of Saint Bonaventure, who introduces the *spirito profetico* (141) himself. And, as Guy Raffa noted, Dante may also be nodding to Joachim's Trinitarian model when mentioning a mysterious third ring in the Sphere of Sun and evoking the Holy Spirit (*Par.* 14.70–78).[21] Steno Vazzana, in one of the few detailed studies of the three *giri* to date, believes that Joachim's diagram was in the forefront of Dante's mind when writing the verses on the Trinity as it was likely that he wanted his final vision to have a mystical image endorsed by a saint's authority.[22]

Another representation of the Trinity in the form of intersecting circles that has been pointed to as a possible source of Dante's vision is a miniature contained in a late thirteenth-century French manuscript (fig. 3). The manuscript, destroyed in a fire in the 1940s, was copied by Adolphe Didron for his book on Christian iconography.[23] While there is no evidence that Dante saw this manuscript, he may have encountered similar illustrations of the Trinity in other treatises, paintings, frescos, stonework, or mosaics. We will return to this manuscript's particular arrangement in the final two sections of this essay.

A third possible source that may have inspired Dante's depiction of the Trinity as three *giri* is a slim, esoteric text that has received little attention by Dante commentators: the anonymous and undated *Liber vigintiquattuor philosophorum* (so designated in the early years of the fourteenth century). The *Liber* exists in three principal redactions and twenty-six manuscript copies, some which may have been reworked by a later thinker. The text is thought by some (Hudry)[24] to have been written in the fourth century, perhaps by a Platonist and not a Christian hand (at least not the earliest version), and by others (Carroll, Poulet, Lucentini)[25] to have been written by a twelfth- or thirteenth-century author. Anna Bagorda's has rigorously explored the potential points of contact between the *Liber* and Dante's *Paradiso*.[26]

In the Middle Ages and Renaissance, the *Liber* was commonly attributed to the mysterious Hermes Trismegistus, also variously to Empedocles, Proclus, Calcidius, Pseudo-Dionysius, Alan of Lille, and the School of Chartres. The *Liber* is only a few pages long. It begins with a single paragraph explaining that twenty-four philosophers (unnamed) gathered together one day and decided that each should devise a sentence-long definition of God. The remaining pages list the results. The first redaction has brief explanatory paragraphs following each sentence. The second has longer explanatory paragraphs, and the third jettisons the explanations altogether, leaving only the twenty-four definitions. The definitions seem to be a summation of Western esoteric philosophy and theology from antiquity through the first centuries of the Christian era.

Out of these twenty-four definitions, the first two circulated most widely in the Middle Ages. Definition I (similar to def. XVII) reads "God is a monad that generates a monad and in itself reflects a flame of love" (Deus est monas monadem gignens, in se unum reflectens ardorem). This definition was conflated with Augustine's triad of *unitas*, *aequalitas*, and Concordia and was most fully theorized in the works of Thierry of Chartres.[27] Definition II (of a decidedly Parmenidean or Empedoclean slant and similar to def. XVIII) says "God is an infinite sphere whose center is everywhere, and circumference is nowhere" (Deus est sphaera infinita cuius centrum est ubique, circumferentia nusquam). Some, like Augustine and Thomas of York, condemned what were thought to be Hermetic (and heretical) ideas; others, such as Abelard and John of Salisbury, saw in it prefigurations of Christ (see Lucentini); and others, such as Boethius, Alan of Lille, Albertus Magnus, and Thomas, explored similar paradoxical modes for describing God—either drawn from the *Liber* or from the texts that inspired it, depending on the date of the work. Definition II in particular has had an active life in Dante commentaries, being (wrongly) attributed to (an unspecified) passage in Augustine's *De civitate Dei* by Pietro di Dante (1559–64), who stated that Augustine ascribed it to Hermes Trismegistus. Benvenuto da Imola (1375–80), Johannis de Serravalle (1416–17), Cristoforo Landino (1481), Niccolò Tommaseo (1837), and others have followed this attribution in their commentaries on Dante's Trinitarian circle/spheres.[28]

Carroll, Poulet, Hudry, Lucentini and Bagorda are among the few scholars to have explored the *Liber* as a possible source for Dante's depiction of the Trinity. In addition to the first two definitions, they have noted other parallels, such as definition XIX: God as "Forever immovable within movement" (Deus est semper movens immobilis: compare *Par.* 24.131; *Par.* 27.109–14); and definition VII: God as "Beginning without beginning, process without change, end without end" (Deus est principium sine principio, processus sine variatione, finis sine fine: compare *Par.* 26.16–18; *Par.* 33.111–14). Dante's exclamation after encountering the Trinity in *Paradiso* 33 could be considered the culmination of this rich series of motifs, which, if not directly inspired by the *Liber*, likely shared its sources: "O luce etterna che sola in te sidi, / sola t'intendi, e da te intelletta / e intendente e te ami e arridi!" (33.124–26). Dante finally grasps what Love cryptically stated to a forlorn, young Dante in the *Vita nuova* (12.11): "Ego tanquam centrum circuli cui simili modo se habent circumferentiae partes: tu autem non sic" (I am like the center of a circle, to which the parts of the circumference have a similar relation: you, however, are not). It is as if in order to "see"—even if only in part—the mystery of the Trinity, Dante had to enter into the paradoxes of the circle; and to "get it," he would have to *become like a circle*, perhaps even a sort of *imitatio Christi* of the Incarnation in the second *giro* (*Par.* 33.130–33)—the form in which every point on its circumference is equidistant from its center (see *Par.* 13.50–51), that is, every part is in perfect relation with its origin, Love.

## Motion

There is an argument to be made that the *giri* of the Trinity are *not* moving, as God is described in both classical philosophy and Christian doctrine as perfectly still: the "unmoved mover" (*Conv.* 2.3), he who "tutto 'l ciel move, / non moto" (*Par.* 24.131–32). But the Pilgrim may have needed to see the *giri* in motion in order to see them, or distinguish between them, at all. As alluded to earlier, the word *giro* could, in fact, imply motion; and it seems even more likely that Dante had motion in mind given his use of the terms *circulazion* (127) and *circunspetta* (129) in the verses that follow.

If we read Dante's *giri* as "roundnesses in motion," what are the possible ways in which roundnesses could move? The poet certainly does not say anything specific about motion; nor does he leave any clues. Each *giro*, if flat like a circle or disc, could be turning around its own central axis like a pinwheel or frisbee; or it could be spinning on an edge like a coin dropped on a table. If the *giri* are spheres, then they could look like tennis balls with topspin. If the *giri* are arranged in a triangular format, they could be shifting positions along a single plane, like leaves of a pinwheel or juggling balls in flight; or they could be rotating through a plane to form what looks like a sphere, like a spinning coin or the vanes of a radiometer viewed perpendicularly to its rotational axis.

Furthermore, the *giri* could also be partaking of multiple types of circular motion. Like the earth that rotates around its core and revolves around the sun (albeit elliptically), the *giri* could be simultaneously spinning like ballerinas on pointe and whirling around one other like a trio of dancers. And then there is the question of direction: are they spinning clockwise or counterclockwise with respect to the Pilgrim's gaze? Are they whirling around him, or is he orbiting them? Is there a centripetal force pulling the *giri* in toward a center? In the spaceless Empyrean (which "è ogne parte là ove sempr' era, / perché non è in loco e non s'impola," *Par.* 22.66–67), perhaps these details are moot. While we cannot glean anything as to the nature of the *giri*'s motion from the verses under consideration, we can infer a few things. Given the Nicene dictate of 381 that specified how the Son proceeded from the Father but should not be thought of as separate from or inferior to the Father, and the Latin version of the Nicene-Constantinople Creed of 381, which clarified that the Holy Spirit proceeded from *both* the Father and the Son (the heart of the *filioque* controversy) and was equal to and unified with them, all three *giri* would have to be spinning and/or revolving at equal speeds and in the same direction—although a case could be made that the Son, as the reflection or *lumen de lumine*, rotates in the direction opposite to that of the Father. Second, since the Pilgrim has fixed his gaze upward (50) at God and is unable to release it (76–78), it seems unlikely that he is orbiting the *giri*, although by the final lines of the poem when the Pilgrim finally grasps the hypostatic union, however, the possibility remains that he has moved close to the God-Point and is now orbiting It with his *disio* and *velle* (143). And thirdly, if by *d'una contenenza* Dante implied that the *giri* are all occupying a single space, it would be impossible for each *giro* to orbit the Pilgrim as he would become their axis, unless they were hovering over his head.

The specific choreographies, single or multiple, that may have passed through Dante's mind with respect to the *giri*, we cannot know. Only a few twentieth-century scholars have posited that Dante imagined the *giri* to be moving in some way, using as support for their claim his choice of a term (*giro*) implying "turning or spinning," and the word *circulazion* in verse 127.[29] Singleton and Freccero point to a long tradition in both classical and Christian thought that holds that motion implies intellection, and that God is pure act, pure intellection, and the origin of all motion, even if he is perfectly still at his core.[30]

Many medieval theologians held that God's thought moved circularly (human thought linearly, and angelic thought in a spiral) and that to contemplate God we too needed to set our mind moving in a circle.[31] One tercet after his description of the *giri*, Dante professes his awe upon recognizing the "circularity" of God's thought, that is, the circular path of his self-knowing and self-loving: "solo in te sidi, / sola t'intendi, e da te intelletta / e intendente te ami e arridi!" (*Par.* 33.124–26). It will not be until the *fulgore* that the Pilgrim will fully grasp this paradox, and the *giri*, if we believe they are moving, could serve as the final training wheels he needs before making the leap to enact the *rota ch'igualmente è mossa* (33.144) by the love that moves the sun and the other stars.

As noted earlier, Dante's use of the Latinism *circunspetta* (33.129) with regard to the path of the Pilgrim's eyes when scanning the *giro* of the Son, may further support a theory of the *giri*'s motion. *Circunspetta* has been translated variously as "contemplated" (Longfellow, 1867), "regarded" (Norton, 1891–92), "surveyed" (Grandgent, 1909–13), "scrutinized" (Singleton, 1975), "watched" (Mandelbaum, 1984), "gazed on" (Hollander, 2007), and "surveyed" (Durling and Martinez, 2011).[32] But *circunspetta* could be interpreted or translated as "carefully looked around and along with," as it implies moving *carefully around* something. If the *giro* of the Son itself is moving, then the Pilgrim's eyes, while moving carefully around to look at it, would also be moving *along with* it. While some commentators have noted that the Latin *circumspicere* indicates a "looking around" in the sense of moving one's eyes around to see an entire thing,[33] none have suggested that the Pilgrim's eyes are moving *along with* the *giri*, perhaps because they would first have to posit that the *giri* themselves were moving.

Yet in their very shapes, circles and spheres evoke both eternal stillness and eternal motion. On the one hand, they are free from a beginning and an end; on the other, they are in themselves both beginning and end, like the *giro* of Christ, the alpha and omega. While the Creed does not speak of the Trinity's motion, and God is often described in classical philosophy and medieval Christianity as perfectly still, mystical treatises spoke of God's thought as moving in a circular path, continually generating the Son and breathing the Holy Spirit. The dynamic words that Dante used in describing the Trinity—*giro*, *circulazion*, and *circumspetta*—hint at the possibility that he could have imagined the *giri* as somehow in motion. If that were the case, at least one quandary—that of how to reconcile one reading of *di tre colori* and *d'una contenenza* (33.117)—could be resolved, as we will see in the "Configuration" section below.

## Size

As the Pilgrim's eyes move around the *giri*, and even perhaps along with them, it would seem that the Trinity has manifested itself to Dante in a size he can grasp—neither too tiny nor too immense. But does it even makes sense to speak of size given the paradoxes of the Trinity, of a God that is both a still point and circumscribing the entire universe, and the Empyrean as a place that is outside of space? Renaissance commentators Antonio Manetti, Alessandro Vellutello, and Galileo earnestly attempted to calculate the size of Dante's Hell, and the Pilgrim himself endeavored to estimate dimensions of many things on his journey: from the height of Lucifer to the diameters of the angelic circles. The *giri* of the Trinity, however, *per forza* evade any such approximations.

All that Dante says regarding the dimension of the *giri* is that they are *d'una contenenza* (33.117). What *contenenza* means in this context, however, is ambiguous. Most commentators have interpreted it to mean "of the same size," that is, that each *giro* has the same diameter length and thus each *contains* the same amount of area. A number of early commentators, however, did not think of the expression in terms of dimension, but rather held it to mean that the *giri* are all three "of the same substance/essence," with the *colori* indicating the different attributes of each person.[34] A few scholars have interpreted the phrase to signify that the *giri* are "together *contained* in one space."[35] While it is clear that Dante wished to recall the unity and equality of the persons of the Trinity (the substance and essence) with *d'una contenenza*, these various interpretations of "containing" and "contained" lead us to imagine quite different arrangements for the *giri*. We will pursue the various ways in which *d'una contenenza* has been and can be envisioned in the "Configuration" section below.

Before leaving the question of *d'una contenenza*, there is a curiosity we would like to consider. In a number of early manuscripts copies of the *Commedia* there is a variation of the word *contenenza* that occurs with some frequency: *contingenza*. Petrocchi has noted this, but rejected it as an "alterazione poligenetica" (an error made by numerous scribes independently); the variant is not discussed by commentators before or after Petrocchi.[36] We do not debate Petrocchi's conclusion, but would like to point out how a copyist could have easily made such a mistake, and why this variation is particularly relevant to our discussion. The word *contingenza* carries within it a number of meanings. The definition of the term that likely first comes to mind is the Aristotelian one that Dante used a number of times in *Paradiso* (it only appears in this canticle), as well as in *De vulgari eloquentia*, *Monarchia,* and *Epistola* 12: "non-necessity, possible occurrence, chance, subject to change, or dependence."[37] Such connotations, however, could not possibly have been what Dante would have associated with the Trinity.

Another use of the word, however, makes far more sense in this context, as it denotes "relationship, affinity, or nature." The Latin *contingere*, similar to *tangere*, signifies a "contact with," that is, "touching, sharing, bordering on," but it also refers to "coloring or imbuing" (*tingere*). In mathematics, "contingency" is equivalent to tangency between two or more objects. Perhaps the early copyists who penned *tre giri d'una contingenza* imagined Dante's Trinity as three roundnesses tangent to one other, comprising a single form, and tinting one another.

## Color

If *contingenza*, with its meaning of touching and tinting together were Dante's intended word, rather than *contenenza*, the question as to how to think about the *giri*'s three colors (117) would be quite a different one. We could imagine each *giro* sharing three colors, perhaps flowing one into the other, with the *giro* of the Holy Spirit the most red of the three, but containing the other two, unspecified rainbow colors. A trefoil knot—common in religious art and architecture of the Middle Ages—would be an excellent format for this co-tinting (fig. 4, in gallery).[38]

On first glance, one might imagine that *tre colori*, following immediately after *tre giri,* indicate that each *giro* was a unique color. Given the specification in verse 119 that the "third" *giro*, the Holy Spirit, seemed like fire—and recalling that the Holy Spirit was associated with charity, love, and the color red in medieval color symbolism[39]—it would seem likely the case that Dante imagined the Holy Spirit's *giro* to be red, and the Father and Son's *giri* to be two other colors, respectively. Giovanni Busnelli, for example, presents the compelling case for the Father's *giro* to be white and the Son's green.[40] He cites Torraca's deduction that they would echo the three theological virtues (also shown in *tre donne in giro* in *Purg.* 29.121–26) and adds that Dante would have followed Thomas and Saint Basil's discussions of John's Apocalypse (the colors surrounding and comprising God) and the rainbow.

But while Dante alludes to redness by saying that the *terzo giro* "parea foco" (119), he notably omits any specific reference to individual colors with respect to the *giri* of the Father and the Son. Instead, he describes an atmospheric phenomenon that most commentators from the earliest to the present have assessed to be the double rainbow, *l'un da l'altro come iri da iri / parea riflesso* (118–19). It appears that Dante wanted to avoid attributing a single color to each *giro*, perhaps in order to avert the criticism that Joachim of Fiore received by giving the three persons three separate colors in his illustration of the Trinity; perhaps to indicate that he could see three colors, but could not quite locate the origin of each; or perhaps because the three *colori* were intended merely to represent the three distinctions between the three persons (the generating, the generated, and love), and the rainbows evoked are symbolic of generation (Father of Son), rather than meant to indicate that they were rainbow colored, or even each colored singularly with one of the rainbow's hues.[41]

Yet even if Dante did not wish to associate himself with colored rings in the Trinity illustration attributed to Joachim in the *Liber figurarum* (that is, if he even saw this illustration), that is not to say he would not find it suggestive. Figure 2 (in the gallery) shows three interlocking rings on a horizontal plane to be "read" from left to right with the Father as the first ring (green), the Son in the center (blue), and the Spirit as the last ring (red). Why green for God and blue for the Son? It is not common in medieval color symbolism to attribute green to God and blue to the Son, the way green was associated with the virtue of hope, white with faith, and red with charity or love. In Revelations 4:2–3—a text from which a number of images in this canto were likely drawn—however, there is a description of an emerald-like rainbow surrounding God, who Himself is likened to jasper (which can be red or yellow-green) and sardine stone. In Joachim's *Expositio* of John's apocalyptic vision, in contrast, he discusses the rainbow, attributing to it three colors, "viridum . . . caeruleum vel aereum . . . rubicundum" (green, blue-purple or sky, and red). Vazzana has proposed that Joachim chose green for the Father because the Father creates and generates, like nature; blue for the Son because the Son descended from the sky; and red for the Holy Spirit to signify the love between Father and Son.[42] Interestingly, however, elsewhere in Joachim's *Expositio* he associates the Father with *jaspidis* (jasper—yellow-green?), the Son with *sardinis* (red), and the Holy Spirit with *smaragdi* (emerald green).[43] What is more, in the same illustrated plate, there is a smaller version of this diagram to its upper right where we see the Holy Spirit's red circle in between the Father and Son and a caption explaining the relationship between the three persons; and in another diagram above that one, a diagram focuses on the Three Ages, with the red ring first, then the blue and then the green. Dante's Trinity, in a sense, allows for all three of these arrangements by *not* indicating how, exactly, the *giri* are oriented. What Dante does do, though, that Joachim's primary illustration does not, is stress how the *terzo giro* seemed to proceed from both the Father's and the Son's *giri*, as co-breathed by each of them in equal measure (120).

Dante has referred to this "co-breathing" before. In *Paradiso* 10, upon entering into the Sphere of the Sun he invites the reader to look up to see the Heavens and contemplate God's triune nature: "Guardando nel suo Figlio con l'Amore / che l'uno e l'altro etternalmente spira" (1–2)—a depiction nearly perfectly repeated in *Paradiso* 33.120, and alluded to in *Paradiso* 13.55–57 ("quella viva luce che sì mea / dal suo lucente, che non si disuna / da lui né da l'amore ch'a lor s'intrea"). The Father and Son eternally breathe together, *con-spire*, to express the Holy Spirit, that is, love for each other and for all of the creation. Thus, while imagining individual colors for each *giro* that could indicate the unique attributes of the three persons of the Trinity—power, knowledge, and love (recalling the words on the Gates of Hell, *Inf.* 3.5–6)—Dante's depiction leans toward representing the oneness of the three, with the Father and Son reflecting one another *come iri da iri* (118), and sharing, one would assume, the same colors. The *giri* of the Father and Son, as Dante describes them, present a double analogy, as they are both in a relationship of reflection and one of co-breathing. The fiery Holy Spirit seems to function as a link between the breaths, and although Dante does not explicitly indicate what is in the space of the reflection between Father and Son, it would be logical for it to be the Holy Spirit as well. Perhaps Dante noticed, in fact, that rainbows have blue on the interior of their bow and red on the exterior. A secondary rainbow flips the spectrum (the Son is a reflection of the Father, *l'un d'altro . . . parea riflesso*) and has red on its inside and blue on its outside. The shared "red" between the bows parallels the shared breathing of *terzo giro che parea foco*, the Holy Spirit.

From the time of the earliest commentators, most have taken the metaphor of reflection of the Father by the Son, described as *come iri da iri*, to indicate a double rainbow. Dante had certainly contemplated such phenomena, as can be seen in his description of a double rainbow in *Paradiso* 12.10–15, which he used as an analogy to the echoing of motion and song between the two rings of saints in the Sphere of Sun as well as to recall the pact between God and Noah after the Flood (Gen. 9:13) that such an event would never happen again. While not the first time Dante mentions rainbows in the *Commedia*,[44] it is his first reference to a double rainbow, and it is located, perhaps not coincidently, in the Heaven of the Sun.

Departing from Aristotle's theory of secondary rainbows—which purported, incorrectly, that the secondary rainbow was further from the viewer (and hence fainter)[45]—Dante, like most natural philosophers of his day held: that a secondary rainbow seen outside a primary one was not farther from the viewer, but a reflection of the inner one (it is in fact a double reflection of sunlight inside individual rain droplets). And although the double rainbow continues to be the most accepted interpretation for the analogy of how the Father's *giro* was reflected in the Son's,[46] many scholars have noted that there are a few obstacles to making it a neat fit for the Trinity:[47] a secondary rainbow is fainter than a primary bow; a secondary bow is larger in diameter than a primary one, as it is concentric to it; and rainbows are almost never seen as complete circles (unless from an airplane or a very high mountain),[48] as we usually see them at an angular radius of about 42 degrees centered at a point directly opposite the sun (the lower the sun on the horizon, the more of the semicircle we see). The first two qualities of double rainbows would lead to unacceptable inequalities between the persons of the Trinity;[49] and the third does not lend itself to circular *giri*. Another issue is the large band of dark sky (Alexander's band) that exists between the two rainbows—Dante's Trinity model does not account for that. Of course, in the spacelessness of the Empyrean, the laws of physics do not have to hold, and both bows can be equally bright, equal in size, form a roundness of some sort and share that roundness with the fiery Holy Spirit *giro*, free of any dark band of space in between them.

André Pézard offered an alternative to the double rainbow analogy: a parhelion,[50] which is similar to a sunbow or halo and occurs when ice crystals in high clouds bend light rays at a minimum deflection of 22 degrees. A rainbow-like ring forms around the sun with two bright bursts, or "sun dogs," that appear on either side of the sun, and sometimes above it. Pézard noted how this model would account for the *iri da iri* reflection without the issue of lack of concentricity, or one bow being fainter than the other, and also how the fiery *giro* of the Holy Spirit links the Father and Son, as the parhelion's color on the side closest to the sun is red (the rainbow's is blue) and forms the edge of the ring around the sun. But what about the fact that parahelia are best seen when the sun is just rising or setting? And what about the fact that there is no horizon in the Empyrean? The idea of the Trinity being associated with a phenomenon that is seen low on a horizon is not one Dante is likely to have held. And the sun dogs themselves do not form *giri*, but rather look like little, arching suns. It is more likely, thus, that Dante intentionally (or not) conflated a number of atmospheric and optical phenomena together in devising his Trinity: parts from double rainbows, and perhaps aspects of parahelia. And perhaps he even was thinking of the full-circle halos (sunbows and moonbows) or coronae.[51]

Although Dante uses *iri* (a *hapax* in the *Commedia*) as the word to describe the Son's *giro* reflection of the Father's, Iris's rainbow colors appear in halos and especially in coronae, and halos and coronae's visible circularity makes them compelling models (or partial inspirations) for the *giri*. In attempting to assess how large the circumferences of the concentric circles of angels seemed from where he stood in the Primum Mobile (*Par.* 28), the Pilgrim compares the distance between what looked like a wildly spinning circle of fire (the Seraphim) and the tiny God-point to that of an *alo* from a source of light (23–24). While he may be referring to an atmospheric halo, which is produced by light reflecting and refracting off of ice crystals around the sun, moon, or strong light source; he may also be thinking of a corona, produced by the *diffraction* of light from the sun or moon by small, uniform water droplets (or ice crystals) of a cloud or foggy glass surface. Both meteorological phenomena have a yellowish-red fire color as their innermost visible ring (the opposite of a rainbow, which has bluish-purple as its innermost color), but halos have a single, more pronounced, thin, virtually colorless ring, while coronae have wide, numerous, colorful rings. Interestingly, in *Paradiso* 10 Dante uses the word *corona* in verse 65 to describe the ring of saints forming around him, comparing it to the ring that forms around the moon when the air is dense with fog. Similarly, in *Paradiso* 28, when comparing the size of the Seraphic circle to that of an *alo* (or corona?), Dante mentions that this *alo* girdles the point of light when vapor "carrying it is most dense," which is actually more a condition of coronae than halos. Furthermore, because "*cotanto*" (22) implies a nearness, and *dipigne* (23) a painting or colorfulness, it seems more likely that he is referring to the corona. Coronae (arguably more so than halos), double rainbows, and parahelia certainly offer abundant properties for a visualization of the Trinity.

Another potential rainbow-color rich natural occurrence that Dante may have encountered and one that he might have mixed into his recipe for envisioning the Trinity is the soap-bubble (fig. 5).[52] Within the perfectly spherical single bubble there is an iridescent swirl of colors accompanying what looks like a double reflection: a reflection of something on the convex, outer part of the bubble, and an inverted, perfectly tangent reflection of that thing on the inside of the transparent bubble. The Nicene statement that the Son is *lumen de lumine* would square well with the way these two reflections interact. The three *giri* could be *d'una contenenza* within a single bubble, with two reflected and joined roundnesses and a pronounced, band (fiery red as the Pilgrim saw it) surrounding the whole and "co-breathed" by the two reflections. Of course, in the Empyrean, the bubble would no longer carry its famed evanescence, but be eternal.

Returning to the question of the color of the *giri*, in Dante's time there was no accepted list of colors in the rainbow, or rather the colors of *Iris* or the *arco*, as it was commonly called.[53] Most natural philosophers followed Aristotle, who attributed three colors to the rainbow (purple, green, and red, see *Meteor.* 3.2). Some classical philosophers and medieval theologians attributed four colors to the rainbow as a parallel to the four elements, four seasons, humors, and directions.[54] As H. D. Austin noted, a few medieval church fathers, such as Isidore of Seville and Bede, even went as far as to think of the rainbow in terms of two colors—blue and red—to signal the Deluge and the Last Judgment.[55] And many works of art that used the rainbow as symbolic of Christ in Majesty focused more on the brightness of the rainbow than on its particular colors.[56] On the other extreme are poets (famously Virgil and Ovid), who celebrated the rainbow's thousands of colors,[57] and today's physicists, who say that the number of colors of the rainbow is indeterminate. The optical theory that gives the rainbow seven color bands did not come about until Newton (although Ptolemy held there to be seven colors in his *Optics*, his theory did find many followers). Some commentators have thought that Dante unwittingly anticipated this view in *Purgatorio* 29.73–78 when likening to the rainbow to the long color streaks left in the air as the seven candle flames of the pageant advanced.[58] Others, however, suspect that Dante meant that each streak contained the three ''Aristotelian'' colors of the rainbow.[59] Perhaps Dante did not wish to specify the colors painting the air by the flames, given the disagreements over the number of colors held in a rainbow. Not being able to determine what colors Dante gave the rainbow, however, does not make contemplating the colors of the Trinity's double rainbow *giro* of Father and Son any more difficult. Instead, it further confirms the image's allusiveness, like the rainbow itself.

### Configuration

Aided by our considerations of the variables of shape, motion, size, and color of the *giri*, we can now explore possible configurations of three roundnesses. Which ones and how many of these passed through Dante's mind is impossible to say, but this question gives us the opportunity to play the combinatorialist, or the kabbalist, permuting the components of the Trinity as a means to contemplating divine ineffability. As commentators and our own analysis have shown, the circle is the most likely roundness Dante had in mind when envisioning his *giri*. As such, our discussion of the *giri*'s possible configurations will focus on circles, although with a few references to spheres and other round figures.

Most commentaries on the *giri*'s configuration have imagined them to be arranged in a planar format. Only a few have theorized how the *giri* could, instead, be arranged using spherical geometry.[60] And although very few scholars have discussed these or any other possible arrangements at length, most would agree that while the Pilgrim may have observed the *giri* as arranged in a planar or spherical format, the true *giri* could, actually, have been entirely other—something the Pilgrim would have been unable to perceive fully.

Let us begin by imagining that the Pilgrim saw the three *giri* arranged in a planar format. Three circles on a two-dimensional plane could have a number of configurations: lined up horizontally or vertically, or arranged rotationally around a central axis. The horizontal arrangement evokes Joachim's diagram, but poses the question of the hierarchy from left to right of Father, Son, and Holy Spirit. The vertical arrangement also evokes a hierarchy (top to bottom). Three circles arranged around an axis yield an equilateral triangle, which also risks having a hierarchy, with a Trinitarian person at the top and two at the bottom, or vice versa. If the triangle were revolving, however, top and bottom would become moot. How linearly arranged circles would be moving around each other in two dimensions (and thus escaping the hierarchy question) is much more difficult to imagine.

Medieval Christian iconography struggled with this issue of hierarchy in depicting the Trinity. How does one visually represent co-equality and consubstantiality without one of the three persons at the "top," or in the "center" (and larger and/or higher), or to the "left" of a sequence? As the Nicene Creed explained, the Father generated the Son (the begotten) and the Holy Spirit proceeded from the Father and Son, but the Son and Holy Spirit were not thought to be subsequent or inferior to the Father. Artists did not concur on how to depict the threeness of the Trinity. Visualizations, thus, show a variety of arrangements of the three persons, and while the Father is often placed above the other two in a vertical sequence, there are also many examples of a horizontal sequence with God the Father to the left of the three (the first read by the eye of someone reading left to right, but on the right from the Trinity's perspective), or in the center and a little larger or higher than the other two, or in a triangular diagram with the word for God on the top of triangle—either at the upper left vertex if the triangle has two vertices on the top (the triangle is pointing downward), or on the single vertex, if only one is on the top (the triangle is pointing upward). The illustration of the Trinity in Alessandro Vellutello's edition of Dante's *Paradiso* is one of the more unusual ones we have seen, as it places the Father to the right of the Son instead of the reverse, the Holy Spirit above them in the form of a dove (this decision is not unusual), and includes a fourth component: Mary, equal in size to the Father and Son and hovering below them (fig, 6).[61]

But returning to the Trinity as three components and arranged in a triangular format. When drawn as a diagram, the equilateral triangle was often used to symbolize the Trinity, as its three individual but equal angles together comprise the whole figure. We can see an example of this arrangement in Peter of Poitiers's 1210 *Compendium Historiae in Genealogia Christi* (fig. 7)—which itself may have been inspired by a diagram by Petrus Alfonsi (fig. 8; cf. fig. 1). The Father is on the upper left vertex, the Holy Spirit on the upper right, and the Son on the lower vertex. All three persons are linked to each other and to a central axis/circle (Deus) in two ways: the outer links, which form the triangle between the vertices, each display the phrase "non est." The inner links, which connect the three persons to the Deus circle, have "est" in them. These statements show the unity and the three-as-one, but also that there is a distinction between the three persons, or rather, a distinction formed by their relations to one another (the Father generates, the Son is begotten, and the Holy Spirit proceeds from both). That is, the Father, Son, and Holy Spirit are all God, and God is all three; but the Father is not the Son, the Son is not the Father, and so on. Dante's *giri* may, in fact, have been inspired by a *scutum fidei* diagram of this sort, and if the nodes were rotating, the question of who is on "top" would become moot.

The *scutum fidei* format becomes an interesting model for Dante's Trinity when we think of the description of *d'una contenenza* as meaning both "all contained in one space" (the links form what looks like a strong, single figure) and each *giro* as "containing the same amount of space." If the Father and Son were at the top vertices (vis-à-vis the Pilgrim's perspective), the *iri da iri* reflection would make sense, as would the co-breathing "down" into the Holy Spirit.

The three persons' circles of the *scutum fidei*, however, while clearly linked to each other, are not themselves touching or linked in the way Joachim's diagram shows them to be. It is not clear, in fact, from Dante's description of the *giri* how he intended them to be linked; but they had to be united/unified in some way. Are the three circles' circumferences not touching at all (like the *scutum fidei* model)? Are they tangent to one another? Are they intertwined with each other? Are they entirely overlapping each other? The option of their not touching at all could only work if they were *d'una contenenza* in the sense that they were contained within a singular space, or were linked together like the persons of the *scutum*. The "tangent circles" option is more suggestive than the "non-touching" option, as three tangent circles on a single plane will satisfy the requirement of all three in contact/relationship with each other; and by forming a tricuspid curve (or deltoid) in their center of contact, the "hierarchy" is inverted (fig. 9), further reinforcing visually the equal relationship between the three persons. And if we believe that Dante may have intended the word *contenenza* to be *contingenza,* there would be a further point in favor of the "tangent circles."

But what if, by *d'una contenenza* Dante intended that the three *giri* inhered in one space. As it is likely that he would have thought the three Persons should be of the same size, thus not concentric to one another, if they occupied the exact same space they would fully overlap each other, like a stack of rings (if circles), or a pile of coins (if discs). How, then, could the Pilgrim distinguish between them (their threeness and their colors) if looking at them face on? Were they transparent? Possibly, but he does not convey this in the verses, and the brightness of this vision would likely make the question of transparency moot. Was he viewing the rings from the side, as then he could see three distinct rings? Unlikely, as then it would be difficult for the Pilgrim to see the *effige* within the *circulazion* of the Son, and it would make less sense to use the "squaring of the circle" analogy to convey his inability (without aid of the *fulgore*) to understand the mystery of the Incarnation. So, if by *d'una contenenza* Dante implied inhering in one space (rather than "all of the same size" or even "all of the same substance"), it is more likely that he meant that they all *were contained* within a single space, either inside a larger circle (which he does not mention, but would evoke Alfonsi's tetragram-Trinity in figure 1), or somehow linked.

Romano Amerio offers a theory as to how three circles could all be contained in one, planar space and be neither contained by a larger circle, nor linked in some way. He imagines the *giri* as concentric circles (the only scholar as far as we know to do so), but with areas adjusted for their diameters in order to equalize their size (fig. 10, in gallery).[62] The Father, closest to the center of the "containing circle" is the smallest, but thickest ring. The Son does not follow, but rather is placed on the outer, thinnest ring. The Holy Spirit is placed in the middle ring in order to account for the co-breathing of the Father and Son. In this model, the Son is meant to reflect the Father quite literally like a secondary rainbow does from the primary one. While Amerio does not show it in his diagram nor mention it in his article, one could deduce that he would agree that the color spectrum of the Son's ring is an inversion of the color spectrum of the Father's.

Giovanni Busnelli, on the other hand, came up with a way of imagining three circles inhering in a single space that did not require concentricity, or a larger circle surrounding all three.[63] In his model, the *giri* appear as three "great circles" (a circle that has a diameter the same as the sphere's, like the earth's equator) equally spaced on the surface of a sphere and meeting only the poles (fig. 11). Although Busnelli does not mention this, the great circles of his model recall the maximum circles (meridians) of the celestial spheres, colures, which has often been mistranscribed as *concoluri* for *concolori* in *Paradiso* 12.11. The word *concolori* itself has, according to Rosetta Migliorini Fissi, often been misinterpreted to mean "of different colors," rather than "of equal color," which she argues is what Dante intended it to mean.[64] Perhaps the *tre colori* of Dante's Trinity is meant to echo this notion of equal color (and colures), even with one *giro* appearing a fiery red to the Pilgrim.

Amerio rejects Busnelli's model, as he does not think the *giri* of Father and Son breathing together to form the Holy Spirit from *two* opposite directions (opposite ways in each hemisphere), fits with the Trinitarian Creed, which implies a single direction of breathing-forth.[65] Yet his objection can be countered, as the Father and Son will always be breathing in the same direction with respect to the Holy Spirit, even if not with respect to our (the Pilgrim's) point of view. What is more, the Pilgrim's point of view would further reinforce the *iri da iri* reflection.

Another possible way of arranging three circles in three dimensions is to place all three in such a way that they touch at a single point, their axis. If moving, they would sweep out a three-dimensional figure called the horn torus (fig. 12), which looks similar to a sphere, but has an infinitely small hole at its center, like a donut missing a hole.

With such an arrangement, the Pilgrim could be viewing the three *giri* alternating past him, like bound pages being flipped continuously (see figure 12 for a representation); or, if we think of the *giri* as spheres, then like three tangent balls revolving around their shared axis such that if the Pilgrim is facing them from the side, he sees one ball in front at any one instant (fig. 31). If he is viewing them face on, they will look, instead, like figure 9.

Seeing *into* a volume of the *profonda e chiara sussistenza* (115) does not necessarily indicate that one is seeing volumes.[66] Our bubble model (fig. 5) is an example of this. Thus even if the three *giri* were spheres such as those in figure 13 (or a figure of an even higher dimension), they still could have been perceived by the Pilgrim as three tangent circles (like fig. 9). The most commonly-imagined arrangement of the *giri*, however, has been of three *intertwined* circles, which as knots, are actually occupying three dimensions, but when seen face on, can seem flat. Medieval European art and architecture offer numerable examples of triple interlocking rings, crescents, branches, leaves, and other figures, both three-dimensional and flat. They are generally arranged in one of two layouts: triangular, either in a continuous triquestra/trefoil knot format (see fig. 4), or a figure with threefold rotational and reflective symmetry (like the leaves of the three-leaf clover, see fig. 9); or as a sort of linear braid (like Joachim's diagram and the *terza rima* itself!). Antonio Rossini has done an excellent study of knot images that Dante could have seen in church mosaics, intarsias, and pavements in Florence, Rome, Ravenna, and elsewhere.[67]

The final section of this study is a topological and combinatorial analysis of the layout of three intertwining circles we find particularly compelling vis-à-vis Dante's description of the Trinity: the round figures arranged in a triangular format with rotational and reflective symmetry. By doing this, we are not *de facto* rejecting the other models of the *giri* that have been proposed. Rather, we wish to explore in depth a configuration we find particularly interesting for what Dante envisioned. Of the many possible "link patterns" in this configuration, we isolate two categories that we see as lending themselves most readily to a Trinitarian model. Given Dante's knowledge of geometry and his attentiveness to Christian doctrine—not to mention the special care he must have taken when describing a mystery so central to his faith—we propose that these two, fascinating *nodi* may have contributed to the Poet's image of God as a syzygy of *sono ed este* (*Par.* 24.141).

### 3-links and the *giri*: The compelling Borromean rings and (3,3)-torus link

When we think of all the possible ways three circles could be linked, as we have done so far in this essay, a kind of vertigo begins to set in. If we narrow our focus to just a single arrangement—that of three, linked circles in a triangular format, it turns out that there is a limit to how many ways they could be configured: there are only sixty-four. These sixty-four patterns of "links," as topologists call them, can be sorted into five categories. Of these five categories, there are two, as we shall see, that provide us with particularly intriguing models for Dante's Trinity: the Borromean rings and the (3,3)-torus link (the "torus" here is *not* related to the horn torus mentioned in the previous section and shown in fig. 12).

To count the ways in which three circles can be depicted mathematically in space, we must solicit the tools of knot theory. A *knot* is an embedding of a simple closed curve (think of it as a single string with its ends glued together) in 3-dimensional Euclidean space, denoted as $R^3$.

The trefoil knot (which we reproduce in the gallery in the far-right of fig. 14; other examples of the trefoil can be seen in fig. 4) is the simplest example of a *nontrivial knot*, as it cannot be "unknotted" to produce something that looks like a circle without cutting and reattaching the ends. The first three examples in figure 14 are depictions of *trivial knots* or *unknots*, which are knots that are not "knotted" at all. Even the figure 8 shape on the far-left of figure 14 is a trivial knot or unknot, since we can untwist the single crossing to produce a depiction that looks like a circle. The circle is sometimes called the *standard unknot*, as it is the simplest depiction of a trivial knot or unknot. Herein and throughout, however, we will call circles either circles or, simply, knots.

What do we call a figure, though, that comprises one or more knots? In mathematics, one or more knots joined together is known as an *n-link*: a finite disjoint union of *n* knots, where by *disjoint* we mean that none of the knots are tangent to each other or intersect each other. For example, any three *non-tangent* or *non-intersecting* circles in $R^3$ can be thought of as a 3-link, as can be seen in the two configurations in figure 15. Here are two configurations of 3-links.

And in figure 16 (found in the gallery) are four examples of three circles in **R**$^3$ that are tangent to each other or intersect each other, and hence are *not* 3-links.

The two examples of 3-links in figure 15 (found in the gallery) are called *trivial* since no two knots (that is, standard unknots here, or circles) in either example are *linked* in a manner such that one winds around another. An example of a 2-link that *does* satisfy this winding property is the *Hopf link* shown in figure 17.

When studying links in **R³**, it is helpful first to consider the *projection* of the link onto the plane **R²**. We may think of a projection of a link as the shadow that is cast onto a wall when a light source (say a flashlight) is directed at the link, which is 3-dimensional. See, for example, the projection onto **R²** of the Hopf link in figure 17.

The Hopf link itself is a 2-link and by definition the two knots (here, two standard unknots or circles) constituting it do *not* intersect, like those in figure 16 do. However, its *projection* onto **R²** is *not* a 2-link, because the "shadow observed on the wall" shows two intersecting circles (and as we will recall, no intersecting or tangent figures form a link). A projection loses critical information about the original link from which it was projected—for example, a projection cannot tell us which knots in the link wind around which other knots in the link. Nonetheless, projections are very useful when counting all possible arrangements of 3-links in a particular configuration.

One question we can ask is how many topologically distinct embedding types (explained in the next paragraph) of three circles in **R³** have the following projection onto the plane **R²** (fig. 19). In answering this, we are not discovering anything new for mathematics, but we hope to show why such a question is useful when thinking about Dante's three *giri*.

Let us recall that *this projection* (fig. 19) is *not* a 3-link since the three circles intersect (the circles in the actual 3-link wind around each other, and hence form a link); in fact, any two circles coincide in exactly two points. We have indicated the figure's six intersection points with arrows. The original 3-link, from which figure 19 is its projection, can have an overcrossing or an undercrossing at each of the six intersection points. That gives a total of $2^6 = 64$ possible depictions of 3-links. However, many of these depictions are what knot theorists call *isotopic*, that is, they are considered to be equivalent. Of all sixty-four 3-links that have the projection in figure 19, there are three symmetries that yield isotopic 3-links:

**Symmetry 1**: rotation by 120 degrees
**Symmetry 2**: reflection about the vertical axis of symmetry
**Symmetry 3**: interchanging all crossings

Any sequence of these three symmetry operations on a 3-link produces a 3-link isotopic to the original one. Knot theorists use the term "up to symmetry'' to designate when two or more links are the same; that is, an *n*-link is the same as another *n*-link up to symmetry (see fig. 20).

It can be shown that there are exactly five topologically distinct embedding types of 3-links with the given projection in figure 19. By *topologically distinct embedding types*, we mean that a link belonging to one type cannot be deformed (such as cutting one of the knots) or reflected (such as performing any of the three symmetry operations above) to produce a link from a different type. To distinguish different links, knot theorists often use a tool called a *numerical link invariant*, which assigns a number to each link that does not change when the invariant is applied to any other link in a particular class. One type of simple invariant is the *number of knots* in a link. Clearly a 3-link cannot be isotopic to any 2-link since the numbers of knots that constitute each are different. To distinguish the five topological distinct embedding types of 3-links having the projection in figure 19, we use a numerical invariant called the *linking number*. The linking number measures the number of times that one knot in a 2-link (for example, see fig. 17) winds around the other knot. This number is always an integer and may be positive or negative depending on the imposed clockwise or counter-clockwise orientations of the two knots that form the 2-link. For the purpose of our present study, we do not need to explore the knots' orientations. Instead, we are interested in which part of which circle crosses over or under another circle. Hence we will consider only the absolute values of the linking numbers.

Let us build up to 3-links with knot-components that are circles. First, let us explore all possible 1-links and 2-links with knot-components that are circles. All 1-links are isotopic to the following topological embedding type (fig. 21, in gallery). For 2-links, we can consider the case in which the two knots are not linked and the case in which they are. These represent all possible topological embedding types of 2-links. Figure 22 (found in the gallery) shows one depiction from each particular type. The depiction on the left is called a *2-component trivial link*, and the one on the right is called a Hopf link, like the one in figures 17 and 20. The 2-component trivial link has linking number equal to 0 since no single knot winds around the other knot (you can see that the white ring is placed atop the black ring). The Hopf link, on the other hand, has one knot winding around the other exactly once, so its linking number equals 1 or -1 depending on the orientation of the two knots. Since we are not assigning clockwise or counter-clockwise orientations to the links we are studying, it is sufficient to consider the absolute value of this Hopf link's linking number. Hence, we say that the Hopf link has linking number equal to 1.

We are now ready to answer the question posed earlier. Namely, how many topologically distinct embedding types of three circles in $\mathbf{R}^3$ have the projection given in figure 19? In that projection, each pair of circles meets in exactly two points. So the linking number of each pair is either 0 or 1 in absolute value. Notice that the three circles in figure 19 have both threefold rotational symmetry and reflection symmetry. These are symmetries 1 and 2 (of the three given earlier). To imagine threefold rotational symmetry or reflection symmetry, think of a three-leaf clover without a stem. If we rotate it 120 degrees or reflect it about the vertical axis of symmetry, then we would see the same clover. So to count topological classes of 3-links, it suffices to count how many *pairs* of circles are linked. Due to the two symmetries in the projection, it is irrelevant to take into account which of the three circles is linked to which. Thus there are four possibilities:

**Case 1**: all three pairs of circles are linked
**Case 2**: exactly two pairs are linked
**Case 3**: exactly one pair is linked
**Case 4**: no pairs of circles are linked

In the first case, the 3-link produced is called a *(3,3)-torus link* and has two possible patterns (fig. 23a, found in the gallery, is one such pattern). In the second case, the 3-link produced is called a *3-component chain* and has three possible patterns (of which fig. 23b is one). In the third case, the 3-link produced is called a *Hopf link with split component* and also has three possible patterns (figure 23c, found in the gallery, is one such pattern). Lastly, in the fourth case, there are two distinct embedding types for the 3-link produced: namely, the *3-component trivial link* and the *Borromean rings*. There is only one pattern, up to symmetry, for each of these particular types (see figures 23d and 23e, found in the gallery, respectively). The number of possible patterns in each type given above is noted, although without proof, by Cromwell, Beltrami, and Rampichini.[68]

We remind the reader that there are a total of sixty-four possible depictions of 3-links that have the projection given in figure 19, since at each intersection point we can choose to have the crossing go under or over. Figure 23 accounts for only five of these 64 depictions. From the 59 remaining depictions, any other depiction belonging to the 3-component trivial link or Borromean rings will be isotopic to the two patterns given in figures 23d and 23e. However, from the remaining depictions that are not isotopic to the bottom row diagrams, not all will be isotopic to one of the top row patterns in figures 23a, b, and c. This is because, as stated earlier, these three embedding types have more than one pattern that has the same properties of their type. For example, the Hopf link with split component has three possible patterns. Figure 23c gives only one of these three. Figure 24 reproduces this pattern on the right but on the left shows of one of the other two patterns in this embedding type (that is, the pattern on the left will also have the property that exactly one pair of circles are linked while the third circle is not technically linked to either of the other two).

Both patterns in figure 24 are in the same class of embedding type (*Hopf link with split component*, fig. 23c), since they share the same linking properties: the black and grey circles are linked together, whereas the white circle is linked to neither the black nor the grey circles. However, the one on the left seems "more linked" than the pattern on the right. If we pick up the right pattern by lifting the white circle, or any circle for that matter, the white circle can separate from the Hopf link (it could slide out, leaving the black and grey attached). However, if we pick up the left pattern by the white circle, or any of the three circles for that matter, the whole 3-link would raise along with it (the white circle would not slide out). The white circle is "woven" in between the black and grey circles (which themselves are Hopf linked): it rests above the grey circle and under the black one.

Figures 23d and 23e are curious in that they share the property of not actually being linked through winding. In figure 23d's case, the three circles are not directly (Hopf) linked at all, and in figure 23e's case, they are "woven" together, forming what is known as a *Brunnian link*. Here is where things begin to "circle back" to Dante.

The simplest example of a Brunnian link consists of three knots, and when these three knots happen to be circles, they are called the *Borromean rings*, in honor of the Renaissance Borromeo family crest symbolizing three intimately linked ruling families (the Borromei, the Visconti, and the Sforza), although not all examples of the actual renditions of the three rings in Borromean contexts reveal the properties that the mathematical Borromean rings do.[69]

A curious feature of the Brunnian link is that if any component knot is removed, the remaining knots slide apart from one another. In the Borromean rings example (figures 23e and 25), this means that cutting any one circle yields a linking number of the two remaining circles equal to 0. Another way to put it is that if we were to remove one circle (by cutting it), the whole figure comes apart. *No* two individual rings are directly (Hopf) linked together, yet the whole figure can be raised as one 3-link if we lift any one circle-component. The common braid also functions this way. This "all for one and one for all" pattern is remarkable, and speaks to the mysterious simultaneity of the Trinity's one-and-threeness. Configured as it is in figure 25, the Father and Son could certainly be imagined to be in a relationship of reflection—two parallel rainbow arcs, in fact—if they are the two *giri* at the top, and the Holy Spirit below, co-breathed by them. The Father could not exist without the Son and the Holy Spirit, the Son without the Father and the Holy Spirit, and the Holy Spirit without the Father and the Son—if one link is broken, the other two disassociate, too. Might Dante have seen a figure like the one from the Chartres manuscript (figure 3, in gallery), or heard about this type of link from travellers to Northern Europe, where the link was more prevalent in the Middle Ages? Might he have played with three loops and figured out its properties? If he did try himself to create the Borromean rings with three loops of some material, he would observe that it cannot be constructed without bending each of the loops (see fig. 26 in the gallery). Three triangularly-shaped knots can be linked without bending, as can other shapes, like ellipses; but circles, no. If Dante knew this, the figure might have lost its appeal as model for the Trinity. Or perhaps, on the other hand, he may have thought of the fact that if such a figure were embedded in the spaceless space of the Empyrean, the Trinitarian circles would not have to bend to be joined.

And maybe, while thinking about the linking of three circles, Dante would have known of or discovered on his own the pattern we now call the (3,3)-torus link (fig. 23a). This is a 3-link in which all three pairs of circles are linked, that is, each circle links with each of the other two circles. It is, consequently, the "most linked" of the five possible embedding types, as the other four types all contain at least one depiction of a pattern in the class that can come apart.

No single circle can slide out if we were to pick up the 3-link up by one circle. Every circle is part of a linked pair. While this sounds ideal for a model of the Trinity, and the circles of this 3-link do not need to bend in order to be linked (Hopf-linked, as can be seen in fig. 28), the (3,3)-torus link also has a drawback: if we were to remove (by cutting) a circle, no matter which one we removed, the remaining two circles would stay linked. Joachim's diagram (figure 2, in gallery), while not in a triangular configuration, also has the property that each circle is Hopf-linked to the other two. They do not form a simple chain, in which each circle is only linked to the one next to it (like figure 23b), but like the (3,3)-torus link, if we remove one circle, the other two circles remained linked.

The (3,3)-torus link and Joachim's diagram lack the ''all-for-one and one-for-all'' property. That said, who would ever think to, or be able to, remove a person of the Trinity? Once the Father *giro* generated the Son *giro*—identical to him although ''reflected'' in some way—their link was only further strengthened by the Holy Spirit *giro* that is the love they express for one another.

## Conclusion

While these powerfully suggestive geometric models are strikingly consonant with Dante's verses and the Nicene Creed, they are but two of many possible models of *giri*—some inspired by biblical passages and theological treatises, others by the natural world, and still others by aesthetic principles of balance and symmetry—as we have seen. We dare not propose the Borromean rings or the (3,3)-torus link, or any other configuration we have discussed, to be that which the Pilgrim saw and the Poet envisioned as the Trinity, yet they do share many points of similarity with Trinitarian doctrine and medieval imagery. Knowing the particular mathematical properties of these links helps us to see what would have made models such as those of Joachim and the Chartrian manuscript so appealing to Dante; and utilizing the tools of topology and combinatorics helps us to imagine more vividly the possible shape, motion(s), size, colors, and configuration of Dante's *giri*. These contemplations increase the sense of awe and wonder in the face of Dante's poetry—simultaneously so precise and so open—that so masterfully beckons us to journey, along with him, into a divine paradox. As Bosco and Reggio wrote, ''in verità, Dio resta nascosto in tutto il canto . . . ma il poeta credette obbligo suo dare una conclusione 'visibile' al racconto tutto 'visibile' del suo viaggio oltremondano.''[70]

*Arielle Saiber*
*Bowdoin College*
*Brunswick, Maine*

*Aba Mbirika*
*University of Wisconsin-Eau Claire*

# NOTES

We owe a debt of gratitude to many who collaborated on this interdisciplinary project. We especially wish to thank the topologists Erica Flapan (Pomona College), Charlie Frohman (University of Iowa), and Maggy Tomova (University of Iowa) for their input on the final section of this essay; Mark Peterson (Mount Holyoke College), Guy Raffa (University of Texas-Austin), and David Albertson (University of Southern California) for their generous feedback on this study as a whole; our research assistant Michael Hannaman, and the numerous mathematics students at Bowdoin College who found themselves counting 3-links along with us.

1. Citations of the *Commedia* are from *La Commedia secondo l'antica vulgata*, ed. Giorgio Petrocchi (Milan: Mondadori, 1966–67).

2. Dante also refers to the Trinity in *Par.* 10.1–6; *Par.* 13.25–27, 52–57, 79–87; *Par.* 14.28–29; *Par.* 15.47; *Par.* 24.139–144; *Par.* 31.28; *VN* 29.3; and *Conv.* 2.5.

3. See especially Giacomo Poletto (1894), Carlo Steiner (1921), Attilio Momigliano (1946–51), and Umberto Bosco and Giovanni Reggio (1979). Commentaries cited here and below that are not included in the notes can be found on the Dartmouth Dante Project website: http://dante.dartmouth.edu/. Many other commentaries and lecturae of *Par.* 33 were consulted for this study, but not directly cited. Among those, the following should be noted: Dante Bianchi, ''Commento metrico al XXXIII canto del *Paradiso*,'' *Giornale dantesco* 37 (1936): 137–47; Peter Dronke, ''L'amore che move il sole e l'altre stelle,'' *Studi medievali* 6 (1965): 389–422; Marguerite Mills Chiarenza, ''The Imageless Vision and Dante's *Paradiso*,'' *Dante Studies* 90 (1972): 77–91; Eugene M. Longen, ''The Grammar of Apotheosis: *Paradiso* XXXIII, 94–99,'' *Dante Studies* 93 (1975): 209–14; Beniamino Andriani, ''Dante's Last Word: The *Comedy* as a *liber coelestis*,'' *Dante Studies* 102 (1984): 1–14; Giuseppe Mazzotta, *Dante's Vision and the Circle of Knowledge* (Princeton: Princeton University Press, 1993); Peter Dronke, ''The Conclusion of Dante's *Commedia*,'' *Italian Studies* 49 (1994): 21–39; Angelo Jacomuzzi, *L'imago al cerchio e altri studi sulla* Divina Commedia (Milan: FrancoAngeli, 1995); Cesare Viviani, ''Il misticismo della parola poetica. Riflessioni sul Canto XXXIII del *Paradiso*,'' *Letture Classensi* 30–31 (2002): 47–51; Christian Moevs, *The Metaphysics of Dante's* Comedy (Oxford: Oxford University Press, 2005); Gino Casagrande, ''Le teofanie di *Paradiso* XXXIII,'' *Studi Danteschi* 74 (2009): 199–24.

4. Mario Fubini, *Due studi danteschi* (Florence: Sansoni, 1951), 83.

5. See Carlo Grabher's commentary (1934–36) to *Par.* 33.115–30.

6. Charles Singleton's commentary (1970–75) to *Par.* 33.116–17.

7. See Bernard McGinn, ''Theologians as Trinitarian Iconographers,'' in *The Mind's Eye. Art and Theological Argument in the Middle Ages*, eds. Jeffrey Hamburger and Anne-Marie Bouché (Princeton: Princeton University Press, 2006), 186–207.

8. For an excellent recent study of the relationship between the *effige* and the circle, see Mark Peterson's ''The Geometry of *Paradise*,'' *The Mathematical Intelligencer* 30, no. 4 (2008): 14–19. See also Arielle Saiber, "Song of the Return: *Paradiso* XXXIII," *Lectura Dantis: Paradiso*, ed. Allen Mandelbaum, Anthony Oldcorn, and Charles Ross (Berkeley: University of California Press, forthcoming).

9. Teodolinda Barolini, *The Undivine Comedy: Detheologizing Dante* (Princeton: Princeton University Press, 1992).

10. Guy Raffa, *Divine Dialectic: Dante's Incarnational Poetry* (Toronto: Toronto University Press, 2000), 116 (emphasis added).

11. Among the few commentators who imagine the *giri* as potentially "spheres" (as well as potentially "circles") are Alessandro Vellutello (1544), H. F. Cary (London: J. Carpenter, 1805–1814), Niccolò Tommaseo (1837), Raffaello Andreoli (1856), Natalino Sapegno (1955–57), Daniele Mattalia (1960), Dorothy Sayers and Barbara Reynolds (London: Penguin, 1949–1962—although in the commentary Reynolds refers to them as "circles"), Umberto Bosco and Giovanni Reggio (1979), and Robin Kirkpatrick (London: Penguin, 2006–2007). See also Giovanni Busnelli, "Dalla luce del cielo della luna alla trina luce dell'Empireo," *Studi Danteschi* 27 (1948): 95–116; and Philip McNair, "Dante's Vision of God: An Exposition of *Paradiso* XXXIII," *Essays in Honour of John Humphreys Whitfield* (London: St. George's Press, 1975). On sphere as metaphor, see also Karsten Harries, "The Infinite Sphere: Comments on the History of a Metaphor," *History of Philosophy* 13 (1975): 5–15.

12. Among the pre-twentieth century commentators who *do* associate the *giri* with circles are Benvenuto da Imola (1375–80), Johannis de Serravalle (1416–17), Cristoforo Landino (1481), Baldassare Lombardi (1791–92), H. W. Longfellow (1867), and Giuseppe Campi (1888–93). See also H. D. Austin, "The Three Rings: *Par.* xxxiii, 116," *Philological Quarterly* 17, no. 4 (1938): 408–11.

13. Giacomo Poletto (1894), John S. Carroll (1904–11), Giovanni Busnelli (1948), Siro Chimenz (1962), and Charles Singleton (1970–75), inasmuch as he also proposes "circling" along with "circle," for example, have directly commented on the polyvalence of the term *giro*.

14. In the *Commedia*, the substantive *giro*—also *giri, giron(e),* and *gironi*—has indicated a circular region (*Inf.* 10.4; *Inf.* 11.30, 39, 42, 49; *Inf.* 13.17; *Inf.*14.5; *Inf.* 16.2; *Inf.* 17.38; *Inf.* 28.50; *Purg.* 12.107; *Purg.* 15.83; *Purg.* 17.80, 83; *Purg.* 18.95; *Purg.* 19.38, 71; *Purg.* 22.2; *Purg.* 23.90; *Par.* 31.67; *Par.* 32.36); an object that encircles (like a belt, wreath, crown, or ring: *Inf.* 31.90; *Purg.* 29.121; *Par.* 12.4; *Par.* 14.74; *Par.* 25.130); a twisting/turning/revolving/wheeling (*Purg.* 9.35; *Par.* 8.20, 26, 35; *Par.* 21.138; *Par.* 23.103; *Par.* 28.15, 139); a celestial sphere (*Purg.* 30.93; *Par.* 2.118, 127; *Par.* 3.76; *Par.* 4.34); a horizon (*Purg.* 1.15); and a passing of time (*Par.* 17.96).

15. Citations from the *Convivio* are from *Opere minori*, vol. 1, pt. 2, ed. Cesare Vasoli and Domenico De Robertis (Milan: Riccardi, 1988).

16. Carroll (1904–11). See also Isidore of Seville's discussion of the invisibility of the Trinity to human eyes, except if it takes a corporeal form (*Etymologiae* 7.1.23).

17. Interestingly, in both *Par.* 30 and *Par.* 33 a nursing infant is invoked right before describing the appearance of the *gronda . . . tonda* in *Par.* 30 and the *giri* in *Par.* 33. In the *gronda . . . tonda*'s case (vv. 82–84), the infant is described as hungrily turning his face to milk upon awakening; but in the *giri*'s case (vv. 106–109), it is the tongue of the infant at the mother's breast that is used as an analogy to how few words the Pilgrim has to describe what he is seeing. Another relevant point of comparison between the *gronda . . . tonda* and the *giri* are the many references to reflection.

18. One the most sustained meditations on circular geometry as a key to understanding God is the work of Nicholas of Cusa. See David Albertson, *Mathematical Theologies: Nicholas of Cusa and the Legacy of Thierry of Chartres* (Oxford: Oxford University Press, 2014), esp. 243–52. A seminal study on the circle in Western symbolism is Georges Poulet, *The Metamorphoses of the Circle* (1961), trans. C. Dawson and E. Coleman (Baltimore: Johns Hopkins University Press, 1966). Useful studies of Dante's use of geometry and geometric symbolism include H. D. Austin, "Number and Geometrical Design in the *Divine Comedy*," *The Personalist* 16 (1935): 310–30; Beniamino Andriani, "La matematica in Dante," *L'Alighieri* 13 (1972): 13–27; Silvio Maracchia, "Dante e la matematica." *Archimede* 31 (1979): 195–208; Robert Levine, "Squaring the Circle: Dante's Solution," *Neuphilologische Mitteilungen* 86 (1985): 280–87; Thomas E. Hart, "Dante and Arithmetic," *The* Divine Comedy *and the Encyclopedia of Arts and Sciences. Acta of the International Dante Symposium, 13–16 November 1983, Hunter College, New York*, ed. Giuseppe Di Scipio and Aldo Scaglione (Amsterdam: John Benjamins, 1988), 95–145; Mark Peterson, "Dante's Physics," *The* Divine Comedy *and the Encyclopedia of Arts* *and Sciences. Acta of the International Dante Symposium, 13–16 November 1983, Hunter College, New York*, ed. Giuseppe Di Scipio and Aldo Scaglione (Amsterdam: John Benjamins, 1988), 163–80; Linda Flosi, *The Geometry of Action in Dante's* Commedia: *Lines, Circles, and Angles* (PhD dissertation, Northwestern University, 1991); Ronald Herzman and Gary Towsley, "Squaring the Circle: *Paradiso* 33 and the Poetics of Geometry," *Traditio* 49 (1994): 95–125; and "*Per misurare lo cerchio* (*Par.* XXXIII.134) and Archimedes' *De mensura circuli*: Some Thoughts on Approximations to the Value of p," *Dante e la scienza*, ed. P. Boyde and V. Russo (Ravenna: Longo, 1995: 265–311).

   19. The *Liber figurarum* is a collection of drawings for Joachim of Fiore's treaties thought to be done by Joachim himself. It was compiled shortly after his death in 1202. See Leone Tondelli, Marjorie Reeves, and Beatrice Hirsch-Reich, ed., *Il libro delle figure dell'Abate Gioachino da Fiore*, 2 vols. (Turin: Società Editrice Internazionale, 1990).

   20. See Tondelli, ed., *Il libro delle Figure*; Marjorie Reeves and Beatrice Hirsch-Reich, *The Figurae of Joachim of Fiore* (Oxford: Clarendon Press, 1972); Marjorie Reeves, "The Third Age: Dante's Debt to Gioacchino da Fiore," *L'età dello spirito e la fine dei tempi in Gioacchino da Fiore e nel gioachimismo medievale: Atti del II Congresso Internazionale di Studi Giochimiti, San Giovanni in Fiore, 6–9 settembre, 1984* (San Giovanni in Fiore: Centro Internazionale di Studi Gioachimiti, 1986): 127–39. See also Francesco Foberti, "Il *Libro delle Figure* di Gioacchino da Fiore," *Sophia* 9 (1941): 332–43; P. Francesco Russo, "Dante e Gioacchino da Fiore," *Dante e l'Italia meridionale. Atti del Congresso Nazionale di Studi Danteschi. Caserta, Benevento, Cassino, Salerno, Napoli 10–16 ottobre, 1965* (Florence: Olschki, 1966): 217–30; Steno Vazzana, "*Parvemi tre giri* (*Par.* XXXII, 116)," *L'Alighieri* 24 (1983): 53–61; Pasquale Iacobone, *Mysterium Trinitatis: Dogma Iconografia nell'Italia Medievale* (Rome: Editrice Pontificia Università Gregoriana, 1997); Marco Rainini, *Disegni dei tempi: Il* Liber Figurarum *e la teologia figurative di Gioacchino da Fiore* (Rome: Viella, 2006); and Antonio Rossini, *Dante: Il nodo ed il volume* (Pisa: Fabrizio Serra, 2011).

   21. See Raffa, *Divine Dialectic*, 138–39.

   22. Vazzana, "*Parvemi tre giri*," 58.

   23. Adolphe Didron, *Christian Iconography; or, The History of Christian Art in the Middle Ages* (1843), trans. E. J. Millington, ed. Margaret Stokes (London, 1891) 2:46.

   24. See Françoise Hudry, ed. and trans., *Liber viginti quattuor philosophorum. Corpus Christianorum* 143A (Turnhout: Brepols, 1997); Hudry, ed., *Le Livre des Vingt-Quartre Philosophes: Résurgence d'un texte du IVᵉ siècle* (Paris: Vrin, 2009).

   25. John Carroll's commentary on the *Commedia* (1904–11); Georges Poulet, *The Metamorphoses of the Circle* (1961), trans. C. Dawson and E. Coleman (Baltimore: Johns Hopkins University Press, 1966); Paolo Lucentini, ed., *Il libro dei* XXIV *filosofi* (Milan: Adelphi, 1999); and Lucentini, "Il *Liber vigintiquattuor philosophorum* nei poemi medievali: Il *Roman de la Rose*, il *Granum Sinapis*, la *Divina Commedia*," *Poetry and Philosophy in the Middle Ages: A Festschrift for Peter Dronke*, ed. John Marenbon (Leiden: Brill, 2001): 131–53.

   26. Anna Bagorda, *Il Paradiso e il* Liber XXIV philosophorum: *L'ente divino ai confini di una metafora* (PhD diss, New York University, 2010).

   27. For Thierry's theory of the mathematical Trinity, see David Albertson, "Achard of St. Victor (d. 1171) and the Eclipse of the Arithmetic Model of the Trinity," *Traditio* 67 (2012): 101–44; and *Mathematical Theologies*, 107–18.

   28. See also P. Pomeo Venturi's 1732 commentary on the *giri*, which focuses on definition I.

   29. See especially Daniele Mattalia (1960) and Pasquini and Quaglio (1982).

   30. In Singleton's commentary he says that the *giri* are "spinning and are thus completely active . . . and spinning in this instance (as everywhere in the *Paradiso*) symbolizes intellection and perfection in complete actualization. All is active in God, nothing is passive." See also John Freccero, "Dante's Pilgrim in a Gyre," *PMLA* 76.3 (1961): 168–81, and "The Final Image: *Paradiso* XXXIII, 144," *MLN* 79, no. 1 (1964): 14–27. See also Stephen Gersh, *From Iamblichus to Eriugena. An Investigation of the Prehistory and Evolution of the Pseudo-Dionysian Tradition* (Leiden: Brill 1978).

   31. See Freccero, ibid.

   32. Charles Eliot Norton (Boston: Houghton Mifflin, 1891–1892), Allen Mandelbaum, (New York: Bantam, 1980–1982), John Ciardi (New York: Norton, 1954–1970), Mark Musa (Bloomington: Indiana University Press, 1996–2004), Antony Esolen (New York: Random House, 2002–2004). Other English translations of *circunspetta* are similar: "mus'd" (Cary, 1805–14), "contemplated with mine eyes (H. F. Tozer, 1904), "dwelt on" (Charles Sinclair, 1939), "appeared so to my scrutiny (John Ciardi, 1954), "looked deep into" (Mark Musa, 1984), "studied" (Anthony Esolen, 2004).

33. The few commentators to have noted that the *circunspetta* indicates eyes moving around to see something are C. H. Grandgent (1909–13), G. A. Scartazzini and G. Vandelli (1929), Carlo Grabher (1934–36), Dino Provenzal (1938), Luigi Pietrobono (1946), Natalino Sapegno (1955–57), Siro Chimenz (1962), and Giovanni Fallani (1965), and Robert Durling and Ronald Martinez (Oxford: Oxford University Press, 1996–2011).

34. For *d'una contenenza* meaning "of the same substance/essence," see Chiose ambrosiane (1355?) Francesco da Buti (1385–95), Johannis de Serravalle (1416–17), Cristoforo Landino (1481), Trifon Gabriele (1525–41), P. Pompeo Venturi (1732), Luigi Bennassuti (1864–68), H. F. Tozer (1901), Francesco Torraca (1905), Carlo Steiner (1921), G. A. Scartazzini and G. Vandelli (1929), and Ernesto Trucchi (1936).

35. For *d'una contenenza* meaning "together contained in one space," see John Carroll (1904–11), Enrico Mestica (1921–22), Charles Singleton (1970–75), Emilio Pasquini and Antonio Quaglio (1982—although they also propose "of the same size"), Mark Musa (1996–2004), Robert Hollander (2000–2007), and Robert Durling and Ronald Martinez (1996–2011).

36. See Petrocchi, ed., *La commedia secondo l'antica vulgata,* vol. 4, 555. He uses as an example Florence, Biblioteca Medicea Laurenziana, fondo principale (Plutei), MS Lau 40 22. From a search of the thirty-three fourteenth and fifteenth-century manuscript copies of the *Commedia* currently available on the Società Dantesca Italiana's website (http://www.danteonline.it/italiano/codici_indice.htm), I have found *contingenza* in the following manuscripts: Florence, Biblioteca Nazionale, Fondo Nazionale II.I.36 / Magl. VII 1032 Strozzi 1281; Florence, Biblioteca Riccardiana, Riccardiano 1010; Florence, Biblioteca Riccardiana, Riccardiano 1094; and Perugia, Biblioteca Comunale Augusta, 240/283. Of the seven, fourteenth-century, key manuscripts included in Prue Shaw's *Dante Alighieri: Commedia; A Digital Edition* (Birmingham: Scholarly Digital Editions; Florence: SISMEL Edizioni, 2010)—two of which are available on the Società Dantesca Italiana's website—there are no occurrences of *contingenza*.

37. Dante uses various forms of the word *contingenza* in the following instances, all implying something that is non-necessary, possible, generated, and/or mutable: *Par.* 13.63–64, 99; *Par.* 17.16, 37, and *Par.* 25.1. See also *DVE* 1.3.1, 9.9; 2.13.9, 14.2; *Mon.* 1.4.2, 5.4, 11.4, 12, 13.7, 15.7; 2.2.3; 3.3.4, 4.4, 6, 11, 15.4; and *Epist.* 12.

38. On the trefoil knot in medieval art and architecture that Dante likely saw, see Antonio Rossini, *Dante: Il nodo ed il volume* (Pisa: Fabrizio Serra, 2011). Mathematically, a trefoil knot is the simplest example of a "nontrivial" knot ("nontrivial" meaning that one cannot untie it without cutting and gluing it back together to get a circular depiction, for instance).

39. On the Holy Spirit's connection to fire, love, and the color red, see Peter Dronke, "Tradition and Innovation in Medieval Western Colour-Imagery," *Eranos-Jarbuch* 41 (1972): 55–107. In the *Commedia*, see especially *Purg.* 6.38; *Purg.* 29.122–23; *Purg.* 30.33; *Par.* 5.1; *Par.* 7.60; *Par.* 21.88; *Par.* 24.145–47; *Par.* 30.52–54; *Par.*31.13.

40. Giovanni Busnelli, *Il concetto e l'ordine del* Paradiso *dantesco* (Florence: Città del Castello, 1911), 259–68.

41. For the idea that the *iri da iri* metaphor implies generation and not the *giri* having rainbow colors, see Vazzana, "*Parvemi tre giri,*" 56. Perhaps not coincidently, in *Purg.* 25.91–96 Statius uses the generation of a rainbow as an analogy to the generation of the soul.

42. Ibid, 58.

43. See Joachim, *Expositio in Apocalypsim,* c. 4 and 21, and Busnelli, 262–63.

44. Other rainbows in the *Commedia* are in *Purg.* 21.50; *Purg.* 24.91–94; *Purg.* 29.73–78; and *Par.* 28.32.

45. For Aristotle's comments on the rainbow, see his *Meteorologica* 3.2–5.

46. See Romano Amerio's argument on how to make concentric circles (like the arcs of the double rainbow) fit with the Creed indicating that they are of the same size, in "Nuove interpretazioni della suprema visione? *Par.* XXXIII, 105 e segg.," *Convivium* 2 (1947): 181–92.

47. See, for example, Rosetta Migliorini Fissi, "*Come iri da iri* (*Par.* XXXIII.118," *Medioevo e Rinascimento* 12.9 (1998): 49–79; and Simone Tarud Bettini, "Dante e il doppio arcobaleno: Una nota su poesia dantesca e scienza aristotelica," *L'Alighieri* 29 (2007): 143–54.

48. Dronke notes that Dante knew that the rainbow was never a perfect circle. See "Tradition and Innovation in Medieval Western Colour-Imagery," 105.

49. On Dante's knowledge of the optics and physics of the rainbow, as well as suggested sources for the use of rainbow imagery in his texts, see especially H. D. Austin, "Dante Notes XI: The Rainbow Colors," *Modern Language Notes* 44.5 (1929): 315–18; Peter Dronke, "Tradition and Innovation in Medieval Western Colour-Imagery" *Eranos-Jahrbuch* 41 (1972): 55–107; Simon Gilson, "Dante's Meteorological Optics: Refraction, Reflection, and the Rainbow," *Italian Studies* 52 (1997): 51–62 and *Medieval Optics and Theories of Light in the Works of Dante* (Lewiston: Edwin Mellen Press, 2000); as well as Migliorini Fissi, "*Come iri da iri* (*Par.* XXXIII.118)," and Bettini, "Dante e il doppio arcobaleno: Una nota su poesia dantesca e scienza aristotelica." For a foundational study on the cultural and scientific history of the rainbow from antiquity through the mid-twentieth century, see Carl Boyer, *The Rainbow: From Myth to Mathematics* (New York: Thomas Yoseloff, 1959), and more recently, Raymond L. Lee, Jr and Alistair B. Fraser, *The Rainbow Bridge: Rainbows in Art, Myth, and Science* (University Park: The Pennsylvania State University Press, 2001).

50. André Pézard, "La vision finale du Paradis," *Mélanges de littérature comparée et de philologie offerts à Mieczyslaw Brahmer*, ed. Charles Vincent Aubrun (Warsaw: Éditions scientifiques de Pologne, 1967): 393–402.

51. For halos (sunbows and moonbows) and coronae, see *Par.* 10.67–69 (arguably a moonbow) and *Par.* 28.23 (an atmospheric halo or corona).

52. People had been making soap with animal fat or vegetable oil and ashes from well before the time of Ancient Rome, although it was used mainly for washing clothing, not personal hygiene, until the late Renaissance. See Hugh Salzberg, *From Caveman to Chemist: Circumstances and Achievements* (Washington, D.C.: American Chemical Society, 1991), 6–7, 74.

53. Dante uses a number of patronymic metaphors in the *Commedia* to indicate the rainbow: "figlia di Taumante" (*Purg.* 21.50), Delia's "cinto" (*Purg.* 29.78), Juno's "ancella" (*Par.*12.12), and "messo di Iuno" (*Par.* 28.32).

54. See Boyer, *The Rainbow: From Myth to Mathematics*, 48 and 74, who cites Empedocles, Democritus, Plato (*Timaeus*), Isidore of Seville (*De natura rerum*), Bede (*De natura rerum*), and Rabanus Maurus (*De arcu coelesti*).

55. Austin, "Dante Notes XI: The Rainbow Colors," 318.

56. See Lee and Fraser, *The Rainbow Bridge*, chap. 2.

57. Virgil, *Aen.* 5.89 and 5.609; and Ovid, *Met.* 6.64. See also Seneca, *Quaestiones*.

58. See Gilson, "Dante's Meteorological Optics" on commentaries that have associated each of the *sette liste* of *Purg.* 29.73–78 with a color of the rainbow. See also Austin, "Dante Notes XI: The Rainbow Colors," 315–16.

59. See Austin, "Dante Notes XI: The Rainbow Colors," 316.

60. See note 12 for commentators who have proposed a spherical arrangement for the *giri*.

61. *Dante con l'espositione di Christophoro Landino, et di Alessandro Vellutello* (Venice: Francesco Marcolini, 1544).

62. See Amerio, "Nuove interpretazioni della suprema visione?," 188.

63. Busnelli, "Dalla luce del cielo della luna alla trina luce dell'empireo," esp. 111–12. See also Francesco Foberti, "Il *Libro delle Figure* di Gioacchino da Fiore," *Sophia* 9 (1941): 332–43.

64. Migliorini Fissi, "*Come iri da iri* (*Par.* XXXIII.118)," 52.

65. Amerio, "Nuove interpretazioni della suprema visione?," 184.

66. Ibid., 186.

67. Antonio Rossini, *Dante: Il nodo ed il volume* (Pisa: Fabrizio Serra, 2011).

68. Peter Cromwell, Elisabetta Beltrami, and Marta Rampichini, "The Borromean Rings," *The Mathematical Intelligencer* 20, no.1 (1998): 55.

69. See Cromwell, Beltrami, and Rampichini's "The Borromean Rings" on the history of this figure and a catalogue of its occurrence in Isola Maggiore.

70. Bosco and Reggio, *Par.* 33.115–32 (1979).

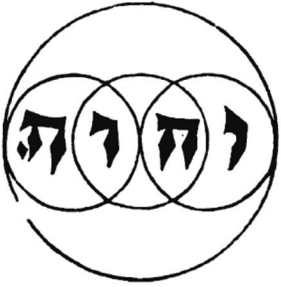

**Figure 1.** Diagram of the Trinity in Petrus Alfonsi's early twelfth-century *Dialogi contra Iudaeos*, ca. 1110 (Tit. 6).

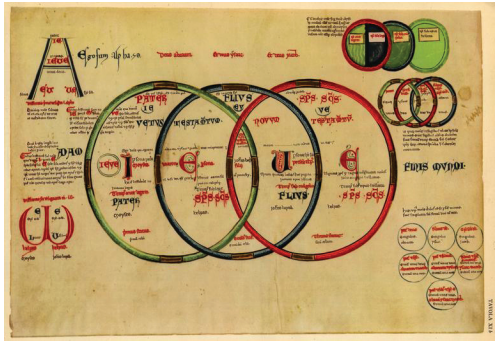

**Figure 2.** Tavola XI, *Liber figurarum*, attributed to Joachim of Fiore. Leone Tondelli, Marjorie Reeves, and Beatrice Hirsch-Reich, ed., *Il libro delle figure dell'Abate Gioachino da Fiore* (Turin: SEI, 1990).

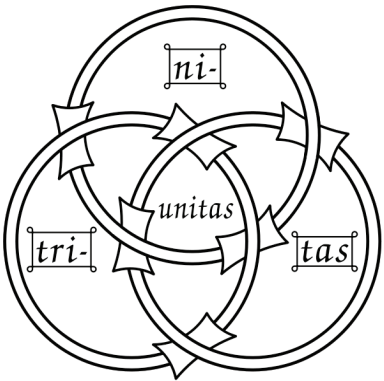

**Figure 3.** A Trinity diagram from a late thirteenth-century French manuscript formerly held at the Municipal Library at Chartres.

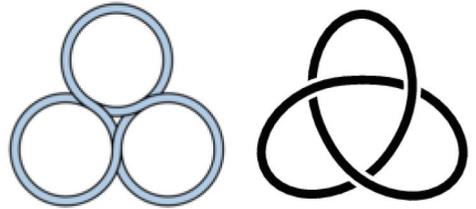

**Figure 4.** Examples of a trefoil knot.

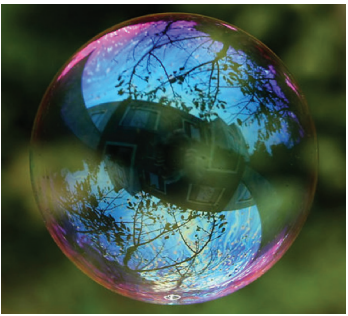

**Figure 5.** Soap bubble. Image from https://sites.google.com/site/thebrockeninglory/.

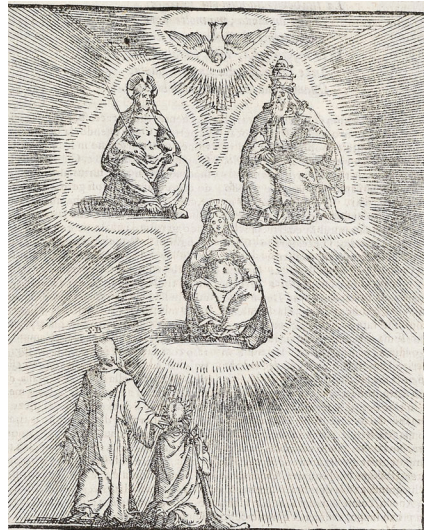

**Figure 6.** *Dante con l'espositione di Christophoro Landino, et di Alessandro Vellutello* (Venice: Marchio Sessa, 1564).

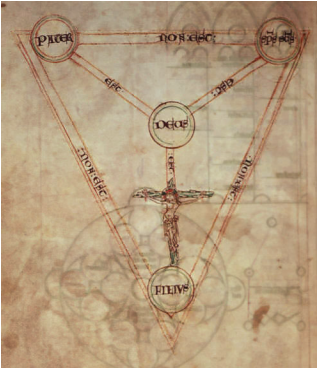

**Figure 7.** Peter of Poitiers, *Compendium Historiae in Genealogia Christi*, ca. 1210. London, British Library, Cotton Faustina MS B.VII, fol. 42v.

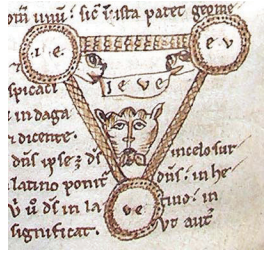

**Figure 8.** Petrus Alfonsi, *Dialogi Contra Iudaeos*, ca. 1110. Cambridge, Cambridge University, St. John's College, MS E.4, fol. 153v.

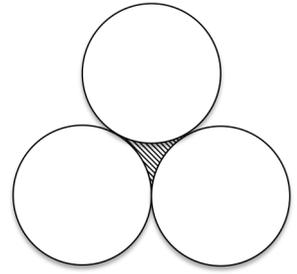

**Figure 9.** Three tangent circles occupying a single plane with tricuspid.

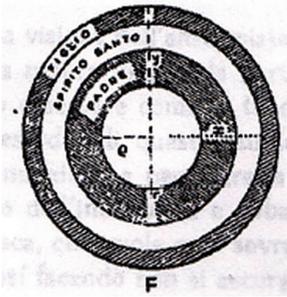

**Figure 10.** Romano Amerio's concentric circles model for the *giri*.

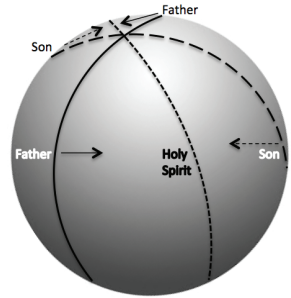

**Figure 11.** The three *giri* as three great circles on the surface of a sphere. Father and Son are breathing forth the Holy Spirit.

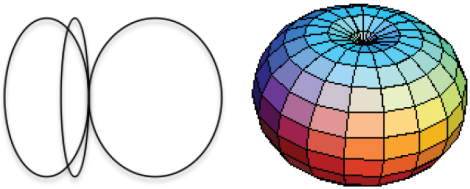

**Figure 12.** *a:* Three circles touching at a single point which, if moving, sweep out a horn torus. *b:* a complete horn torus. From http://mathworld.wolfram.com/HornTorus.html.

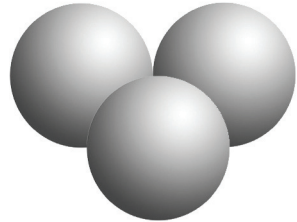

**Figure 13.** Three tangent spheres.

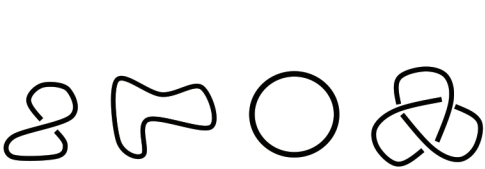

**Figure 14.** Four examples of knots.

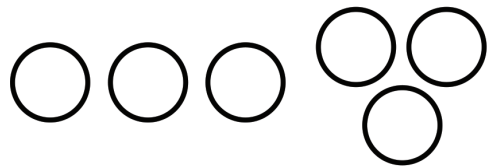

**Figure 15.** Two examples of trivial 3-links.

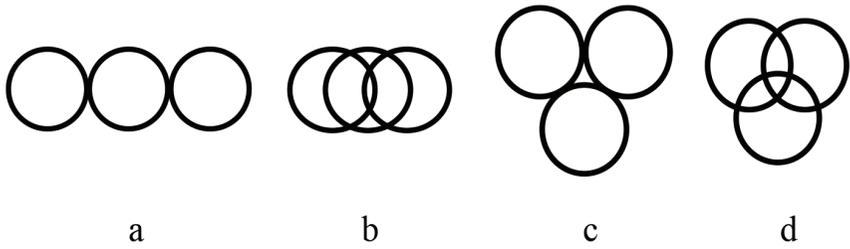

a b c d

**Figure 16.** Four examples of three tangent or intersecting circles that fail to form 3-links.

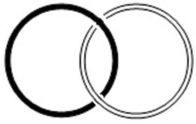 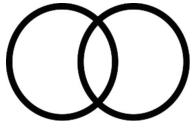 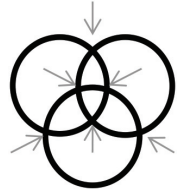

**Figure 17.** A Hopf link.

**Figure 18.** Projection of the Hopf link onto **R**².

**Figure 19.** A projection of a 3-link onto **R**². The projection itself is *not* a 3-link.

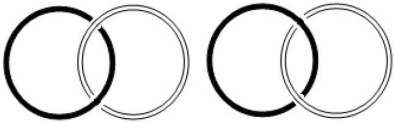 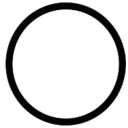 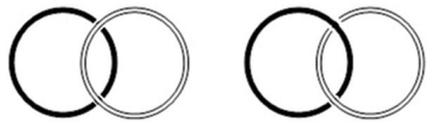

**Figure 20.** Two isotopic 2-links—isotopic because one is merely 180 degrees rotation of the other; as such, they are equivalent.

**Figure 21.** The only embedding type of 1-links with circle-components.

**Figure 22.** The two embedding types of 2-links with circle-components.

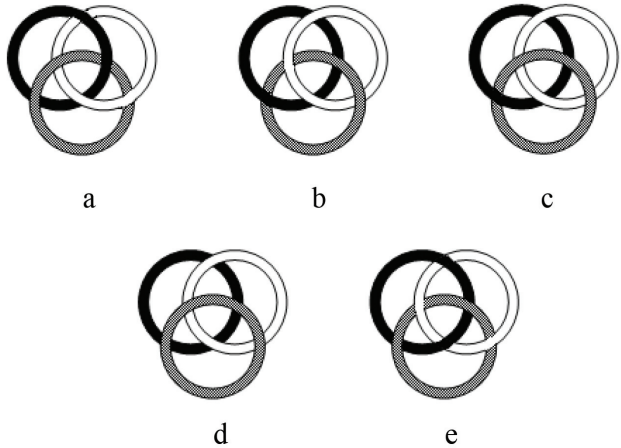

**Figure 23.** Depictions of the five embedding types of 3-links with circle-components. *a:* (3,3)-torus link (1 of 2 possible patterns); *b:* 3-component chain (1 of 3 possible patterns); *c:* Hopf link with split component (1 of 3 possible patterns); *d:* 3-component trivial link; *e:* a Brunnian link, also known as the Borromean rings.

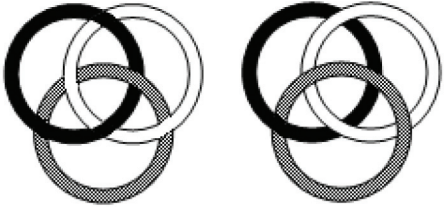

**Figure 24.** Two different patterns of the same embedding type.

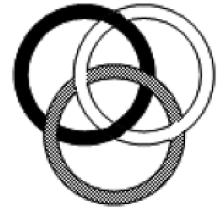

**Figure 25.** The 3-component Brunnian link known as Borromean rings.

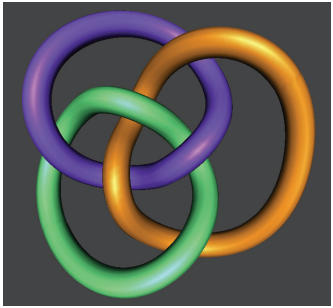

**Figure 26.** A three-dimensional model of Borromean rings, where we can see that the circles have to bend in order to be able to be woven together. Drawn with KnotPlot (http://www.knotplot.com/).

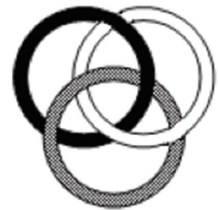

**Figure 27.** The 3-component link known as the (3,3)-torus link. Every pair of circles is linked, making this the most fully linked pattern of the five categories of 3-links in a triangular configuration.

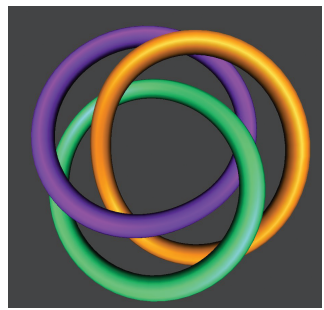

**Figure 28.** A (3,3)-torus link viewed in three-dimensions. The three circles do not need to bend to form this 3-link. Drawn with KnotPlot (http://www.knotplot.com/).